\documentclass{amsart}
\usepackage{amsmath,amssymb,latexsym}
\usepackage{times}
\usepackage[T1]{fontenc}

\font\goth=eusm10
\newcommand\F{\mathcal F}
\newcommand\E{\mathcal E}
\newcommand\Ii{\hbox{\goth I}}

\newcommand\ZZ{\mathbb{Z}}
\newcommand\CC{\mathbb{C}}

\newcommand\N{\mathcal N}
\newcommand\T{\mathcal T}
\newcommand\Oc{\hbox{\goth O}}

\newcommand\Pp{\mathbb P}
\newcommand\PR{\mathbb{P}^r}
\newcommand\Pt{\mathbb{P}^3}
\newcommand\Gr{\mathbb{G}}
\newcommand\Ff{\mathbb{F}}
\newtheorem{theorem}{Theorem}[section]
\newtheorem{coroll}{Corollary}[section]
\newtheorem{lemma}{Lemma}[section]
\newtheorem{proposition}{Proposition}[section]
\newtheorem{definition}{Definition}[section]
\newtheorem{remark}{Remark}[section]

\begin{document} 

\title{Moduli of nodal curves on smooth surfaces of general type}
\author{F. Flamini}

\email{flamini@matrm3.mat.uniroma3.it}
\curraddr{Dipartimento di Matematica, Universita' degli Studi di Roma
- "Roma Tre"\\ Largo San 
Leonardo Murialdo, 1 - 00146 Roma \\Italy}

\thanks{The author is a member of GNSAGA-INdAM}

\begin{abstract}{In this paper we focus on the problem of computing the 
{\em number of moduli} of the so 
called {\em Severi varieties} 
(denoted by $V_{| D |, \delta}$), which parametrize universal 
families of 
irreducible, $\delta$-nodal curves  
in a complete linear system $|D|$, on a 
smooth projective surface $S$ of general type. We determine geometrical and 
numerical conditions on $D$ and numerical conditions on $\delta$ ensuring 
that such number coincides with 
$dim(V_{| D |, \delta})$. As related facts, we also determine 
some sharp results 
concerning the geometry of some Severi varieties.} 
\end{abstract}

\maketitle

\section*{Introduction}\label{S:intro}

Let $S$ be a smooth, projective surface and let $|D|$ denote a 
complete linear system 
on $S$, whose general element is assumed to be a smooth, irreducible curve. 
By the hypothesis on its 
general element, it makes sense to 
consider the subscheme of $|D|$ which parametrizes a universal family of irreducible 
curves having only $\delta$ nodes 
as singular points. Such a subscheme is functorially defined, 
locally closed in $|D|$ (see \cite{W} for $S = {\Pp}^2$ but 
the proof extends to any $S$) and denoted by 
$V_{| D |, \delta}$. It is usually called the 
{\em Severi variety} of irreducible $\delta$-nodal curves in $|D|$, 
since Severi was the first who studied some properties of families of 
plane curves of given degree and given geometric genus (see \cite{Sev}). 

One can be interested in studying the moduli behaviour 
of the elements that a Severi variety parametrizes. This means to understand how the
natural functorial morphism  
$$\pi_{| D |, \; \delta}: 
V_{| D |, \delta} \longrightarrow {\mathcal M}_g$$behaves, 
for each $\delta \geq 0$, where 
$g=p_a(D) - \delta$, $p_a(D)$ the arithmetic genus of $D$ and 
${\mathcal M}_g$ the moduli space of smooth curves of (geometric) genus $g$; 
precisely, the problem is to determine the dimension of the image of 
$\pi_{| D |, \; \delta}$. 

In \cite{S}, Sernesi considered the case $S= {\Pp}^2$. Denote by 
$$ \pi_{n, \, \delta} : V_{n, \delta} \to {\mathcal M}_g$$the functorial 
morphism from the Severi variety of plane irreducible and 
$\delta$-nodal curves of degree $n$ to the moduli space 
of smooth curves of genus $g= \frac{(n-1)(n-2)}{2} - \delta$. Recall 
that $ V_{n, \delta}$ is irreducible (see \cite{H}).

\begin{definition}\label{def:2.expmod}(see \cite{S}) 
The {\em number of moduli} of $V_{n, \delta}$ 
is $dim(\pi_{n, \, \delta}(V_{n, \delta}))$. 
$V_{n, \delta}$ is said to have the {\em expected number of moduli} 
if such dimension equals 
$$min(3g-3, \; 3g - 3 + \rho(g,2,n)),$$where $\rho(g,2,n)$ is the {\em 
Brill-Noether number}.
\end{definition}

\noindent
Of course, when $ \rho(g,2,n)\geq 0$, 
$V_{n, \delta}$ has the expected number of moduli 
$3g-3 = dim({\mathcal M}_g)$ 
when every sufficiently general curve of genus $g$ belongs to it; in such a  
case, this family of plane curves has {\em general moduli}. 
When $\rho(g,2,n)<0$, 
the family $V_{n, \delta}$ does not have general moduli, i.e. it has 
{\em special moduli} and the number $ - \rho(g,2,n)$ determines the 
expected codimension of $\pi_{n,\delta}(V_{n, \delta})$ in ${\mathcal M}_g$. 

\noindent
With this set-up, Sernesi proved the following result:

\begin{theorem}\label{thm:2.mod.sern}
For all $n$, $g$ such that 
$$ n \geq 5 \;\; {\rm and} \;\; n-2 \leq g \leq \frac{(n-1)(n-2)}{2}, $$ $V_{n, \delta}$ has 
the expected number of moduli.
\end{theorem}

\begin{remark}\label{rem:0.1}
\normalfont{
Since $ 3g - 3 + \rho(g,2,n) = 3n + g -9 = 
dim(V_{n, \delta}) - dim (Aut({\Pp}^2)) $, when $\rho(g,2,n) <0$ the fact that 
$V_{n, \delta}$ exactly has the expected number of 
moduli means that its general point parametrizes a curve $X$ which is birationally - but 
not projectively - equivalent to finitely many curves of the 
family, i.e. the normalization $C$ of $X$ has only finitely 
many linear systems of degree $n$ and dimension $2$. 
}
\end{remark}

In this paper, we are interested in the case 
of $S$ a smooth, projective surface of general type. 
In such a case, the expected number of moduli equals 
$dim(V_{| D |, \delta})$ (see Definition \ref{def:5.expmod}).

We determine some general conditions  on $D$, $\delta$ and, sometimes, on the geometry of 
$S$ guaranteeing that such expected number of moduli is achieved (see Theorems 
\ref{thm:5.vannodi}, 
\ref{thm:5.van}, \ref{thm:main}, \ref{thm:m=3,4} and \ref{thm:finale}). 
As a particular case of our more general results, we get 
the following: 
  
\noindent
{\bf Proposition.} {\em Let $S \subset \PR$ be a smooth, non-degenerate complete intersection 
of general type whose canonical divisor is $K_S \sim \alpha H$, 
where $H$ denotes its hyperplane section, $\alpha$ a positive integer and $\sim$ 
the linear equivalence of divisors on $S$. Let $m$ be a positive integer and let 
$X \sim m H$ be an irreducible curve, with only $\delta$ nodes as singular points, of 
geometric genus $g = p_a(X) - \delta$, $\delta \geq 0$. Suppose that $[X]$ is a regular 
point of the Severi variety 
$V_{| mH |, \delta}$ (in the sense of Definition \ref{def:2.reg}).

\noindent
Assume that:  

\noindent
(1) $\delta \leq {\rm dim} (|mH|)$ if 
\begin{itemize}
\item[a)] $\alpha \geq 2$, $m \geq \alpha + 6$, $\delta \geq 1$ or
\item[b)] $\alpha \geq 1$, $m \geq \alpha + 6$, $\delta= 0$; 
\end{itemize}

\noindent
(2) $\delta < \frac{m(m-4)}{4} deg(S)$ if $\alpha \geq 1$ and $5 \leq m \leq \alpha + 5$;

\noindent
(3) \begin{itemize}
\item[a)] $\delta < \frac{m(m-2)}{4} deg(S)$ if $\alpha \geq 2$ and $m = 3, \; 4$ or
\item[b)] $\delta = 0$ if $\alpha =1$ and $m = 3, \; 4$; 
\end{itemize}

\noindent
(4) $\delta <  deg(S) (2 + \alpha) + 
\frac{(r+3-4deg(S))}{2} \chi + \frac{(r-1)}{2} \chi^2$ 
if $\alpha \geq 1$ and $m =2$, where 
$\chi$  is a non-negative integer in 
$[\frac{2deg(S)-1}{r-1}-1 , \; \frac{2deg(S)-1}{r-1})$; 

\noindent
(5) $\delta < \frac{deg(S)}{2} (1+ \alpha) + 
\frac{(r- 2deg(S) + 2)}{2} \chi + \frac{(r-2)}{2} \chi^2$ if 
$\alpha \geq 1$ and $m=1$, 
where  $\chi$ is a non-negative integer in $[\frac{deg(S)-1}{r-2}-1 , 
\; \frac{deg(S)-1}{r-2})$;

\vskip 0,2cm

\noindent
Then the morphism $$\pi_{ |mH|, \delta} : V_{|mH|, \delta} 
\to {\mathcal M}_g $$has 
injective differential at $[X]$. In particular, it has 
finite fibres on each generically regular component of $V_{|mH|, \delta}$, 
so each such component parametrizes a family 
having the expected number of moduli.} 

\vskip 0,2cm

In particular, we have the following:

\noindent
{\bf Corollary.} {\em Let $S \subset \Pt$ be a smooth surface of degree $d \geq 5$ and let 
$[X] \in V_{|mH|, \delta}$ be a regular point.

\noindent
Assume that:  

\noindent
(1) $\delta \leq {\rm dim} (|mH|)$ if 
\begin{itemize}
\item[a)] $d \geq 6$, $m \geq d+2$, $\delta \geq 1$ or
\item[b)] $d \geq 5$, $m \geq d+2$, $\delta= 0$; 
\end{itemize}

\noindent
(2) $\delta < \frac{m(m-4)}{4} d$ if $d \geq 5 $ and $5 \leq m \leq d+1$;

\noindent
(3) \begin{itemize}
\item[a)] $\delta < \frac{m(m-2)}{4} d$ if $d \geq 6 $ and $m = 3, \; 4$ or
\item[b)] $\delta = 0$ if $d = 5 $ and $m = 3, \; 4$; 
\end{itemize}

\noindent
(4) $\delta <  d - 2 $ if $d \geq 5$ and $m =2$; 

\noindent
(5) $\delta < d - 3 $ if $ d \geq 5$ and $m=1$.

\vskip 0,2cm

\noindent
Then the morphism $$\pi_{ |mH|, \delta} : V_{|mH|, \delta} 
\to {\mathcal M}_g $$has 
injective differential at $[X]$. In particular, it has 
finite fibres on each generically regular component of $V_{|mH|, \delta}$, 
so each such component parametrizes a family 
having the expected number of moduli.} 

\vskip 0,2cm

The paper consists of seven sections.
In Section \ref{S:1}, we recall some terminology and notation. Section 
\ref{S:2} contains fundamental definitions and technical details 
which are used for our proofs. Section \ref{S:3} contains the main results of the paper 
(Theorems \ref{thm:5.vannodi}, \ref{thm:5.van}). In Section \ref{S:4} we consider a 
fundamental proposition, which is the key point to determine the results of Sections 
\ref{S:5} and \ref{S:6}. Such theorems focus on cases to which the results of Section 
\ref{S:3} cannot apply. For simplicity, in Section \ref{S:7} we sum up our results in the 
particular cases of Severi varieties of the form $V_{|mH|, \delta}$ on smooth 
complete intersection surfaces of general type or on smooth surfaces in $\Pt$ of degree 
$d \geq 5$.

{\it Acknowledgments:} Part of this paper is contained in 
my Ph.D. thesis, defended on January 2000 at the Doctoral Consortium of 
Universities of Rome "La Sapienza" and "Roma Tre". My special 
thanks go to my advisor E. Sernesi for his constant guide and for having 
introduced me in such an interesting research area. I am greatful to F. Catanese, for 
having pointed out some crucial examples on the subject, and to 
C. Ciliberto, for some  fundamental 
remarks related to the proof of Proposition \ref{prop:chiave}. 
I am indebted to L. Chiantini, L. Ein, A. F. Lopez and 
A. Verra for fundamental discussions which allowed me to focus 
on key examples. I would also like to thank the referee who suggested to use 
Lemma \ref{lem:referee} and projective bundle arguments in order to improve condition (ii) 
in Theorem \ref{thm:5.vannodi} with ''$L$ nef divisor'' instead of the previous 
''$L$ ample divisor''.

\section{Notation and Preliminaries}\label{S:1}
We work in the category of $\CC$-schemes. $Y$ is a \emph{$m$-
fold} if it is a reduced, irreducible and non-singular 
scheme of finite type over $\CC$ and of
dimension $m$. If $m=1$, then $Y$ is a (smooth) \emph{curve}; 
$m=2$ is the case of a (non-singular) \emph{surface}. 
If $Z$ is a closed subscheme of a scheme $Y$, $\Ii_{Z/Y}$ (or $\Ii_Z$) 
denotes the \emph{ideal sheaf} of $Z$ in $Y$ whereas ${\N}_{Z/Y}$ is 
the {\em normal sheaf} of $Z$ in $Y$. When $Y$ is a smooth variety, $K_Y$ denotes a canonical divisor whereas 
${\T}_Y$ denotes its tangent bundle.

Let $Y$ be a $m$-fold and let $\E$ be a rank $r$ vector bundle on 
$Y$; $c_i(\E)$ denotes the \emph{$i^{th}$-Chern class} of $\E$, 
$1 \leq i \leq r$. The symbol $\sim$ will always 
denote linear equivalence of divisors on $Y$. As usual, $h^i(Y, \; 
-):=\text{dim} \; H^i(Y, \; -)$.

If $D$ is a reduced curve, $p_a(D)=h^1(\Oc_D)$ 
denotes its \emph{arithmetic genus}, 
whereas $g(D)= p_g(D)$ denotes its \emph{geometric genus}, the
arithmetic genus of its normalization. For a smooth curve $D$, $\omega_D$ 
denotes its canonical sheaf, i.e. $\omega_D \cong \Oc_D(K_D)$.

\begin{definition}\label{def:2}
Let $S$ be a smooth, projective surface and 
$Div(S)$ be the set of divisors on $S$. An element $B \in 
Div(S)$ is said to be {\em nef}, if $B \cdot D \geq 0$ for each 
irreducible curve $D$ on $S$ (where $\cdot$ denotes the intersection form 
on $S$; in the sequel we will omit $\cdot$). A nef divisor $B$ is 
said to be {\em big} if $B^2>0$. 
\end{definition}

\begin{remark}
\normalfont{
We recall that, given a smooth surface $S$, $N(S)^+$ is the set of divisor classes with 
positive intersection numbers with itself and with an ample class. 
By Kleiman's criterion (see, for example, \cite{H1}), a nef divisor 
$B$ is in the closure of $N(S)^+$. 
}
\end{remark} 

\begin{definition}\label{def:3}
\normalfont{Let $S \subset \PR$ be a smooth surface, $H$ its hyperplane 
section and $D \in Div(S)$. We denote by 
$ \nu (D, H)$ the {\it Hodge number} of $D$ and $H$,
$$ 
\nu(D,H):= (DH)^2 - D^2H^2.$$
By the Index Theorem 
(see, for example, 
\cite{BPV} or \cite{Fr}) this is non-negative since $H$ is a very 
ample divisor.}
\end{definition}

\begin{definition}\label{def:5}
\normalfont{Let $S$ be a smooth, projective
surface. A rank 2 vector bundle $\E$ on $S$ 
is said to be {\em Bogomolov-unstable} if 
there exist $M, \; B \in Div(S)$ and a 0-dimensional scheme $Z$ 
(possibly empty) fitting in the exact sequence
\begin{equation}
0 \to \Oc_S(M)\to \E \to \Ii_{Z}(B) \to 0 \label{eq:2.3}
\end{equation}such that $(M-B) \in N(S)^+$. 
}
\end{definition}

\begin{remark}\label{rem:a}
\normalfont{Recall that $\E$ is Bogomolov-
unstable when $c_1(\E)^2 - 4 c_2(\E) >0$ (see \cite{Bog} or \cite{Reid}).} 
\end{remark}

It is also useful to remind some standard terminology and techniques on 
Severi varieties. 
Consider $S$ a smooth, projective surface and assume that, 
for given $D \in Div(S) $ and $\delta$ positive integer, 
$ V_{| D |, \; \delta} \neq \emptyset$. If $[X] \in 
V_{| D |, \; \delta}$, $N$ will always denote the scheme of 
nodes of $X$, which is a closed zero-dimensional subscheme of 
$S$ of degree $\delta$. From now on, 
denote by 
\begin{equation}\label{eq:1bis}
\varphi: C \to X \subset S
\end{equation}the normalization map of $X$. 
Thus, on $C$ we have the exact sequence of vector bundles 
\begin{equation}\label{eq:5.normalmap}
0 \to {\T}_C \to \varphi^*({\T}_S) \to {\N}_{\varphi} \to 0,
\end{equation}where ${\N}_{\varphi}$ is the {\em normal bundle} of 
$\varphi$. Observe that, if $\tilde{N}$ denotes the pull-back of $N$ to $C$, 
by (\ref{eq:5.normalmap}) we get that 
${\N}_{\varphi} \cong \Oc_C(\varphi^*(D) - \tilde{N})$, so we have  
\begin{equation}\label{eq:2bis}
H^i({\N}_{\varphi}) \cong H^i(\Oc_C(\varphi^*(D) - \tilde{N})), \; i \geq
0.
\end{equation}From 
Horikawa's theory (see \cite{Ho}), 
$H^0({\N}_{\varphi} )$ parametrizes all first-order 
equisingular deformations of $X$ in $S$. Therefore, one gets
\begin{equation}\label{eq:2.Y3}
T_{[X]}(V_{| D |, \; \delta}) \cong H^0(\Ii_{N/S}(D))/<X>,
\end{equation}which is the subspace of 
$H^0({\N}_{\varphi})$ contained in 
$T_{[X]}(| D |)\cong H^0(\Oc_S(D))/<X>$.
\begin{remark}\label{rem:2.regnonreg}
\normalfont{
When $S$ is assumed to be a regular surface, then 
$$(H^0(S, \Ii_{N/S}(D))/<X>) \; \cong  \; H^0(C, \Oc_C(\varphi^*(D) -
\tilde{N})) ,$$which means that all first-order 
equisingular deformations of $X$ in $S$ are in $|D|$, i.e. 
\begin{equation}\label{eq:2.Y2}
T_{[X]}(V_{| D |, \; \delta}) \cong H^0(C, \; {\N}_{\varphi})
\end{equation}
}
\end{remark}

\begin{remark}\label{rem:1.4}
\normalfont{From the exact sequence 
\begin{equation}\label{eq:2.Y5}
0 \to \Ii_{N/S}(D) \to \Oc_S(D) \to \Oc_N(D) \to 0 
\end{equation}and from (\ref{eq:2.Y3}), we get 
$$dim(T_{[X]}(V_{| D |, \; \delta})) 
\geq h^0(\Oc_S(D)) - \delta - 1 = dim(|D|) - \delta;$$the above 
inequality is an equality if and only if the surjection 
$H^1(\Ii_{N/S}(D)) \to H^1(\Oc_S(D))$ is an isomorphism, i.e. if and 
only if $N$ imposes independent conditions to the linear 
system $|D|$. In such a case 
$V_{| D |, \; \delta}$ is smooth at $[X]$ of codimension $\delta$ 
in $|D|$. 
}
\end{remark}We recall the following:
\begin{definition}\label{def:2.reg}$V_{| D |, \; \delta}$ is said to 
be {\em regular} at the point $[X]$ if it is smooth at $[X]$ of dimension 
$dim (|D|) - \delta$. 
Otherwise, the component of $V_{| D |, \; \delta}$ containing 
$[X]$ is said to be a {\em superabundant component}. A component 
of a Severi variety is said to be {\em regular} if it is regular at each 
point, {\em generically regular} if it is regular at its general point. 
\end{definition}We recall that the regularity is a very strong condition, 
indeed it implies that the nodes of $X$ can be independently smoothed in 
$S$ (see, for example, \cite{CS} and \cite{S}).

\section{Basic definitions and technical tools}\label{S:2}

In this section we introduce fundamental definitions 
and remarks which are used to compute the number of moduli of some Severi varieties. 

From now on, $S$ will denote a smooth, 
projective surface of general type, 
unless otherwise specified. Let $|D|$ be 
a complete linear system on $S$, whose general element is supposed 
to be a smooth, irreducible curve. Denote 
by $X$ an irreducible curve in $|D|$ having only 
$\delta \geq 0$ nodes as singularities. As in 
(\ref{eq:1bis}), the map $ \varphi : C \to X \subset S$ 
denotes its normalization, where 
$C$ is a smooth curve of geometric genus $g = p_a(D) 
- \delta$. 

\noindent
We shall always assume that 
$g \geq 2$, for each $\delta \geq 0$. This assumption is not so restrictive for the problems 
we are interested in.

\noindent
With this setup, for each $\delta \geq 0$ one can consider the morphisms:
\begin{equation}\label{eq:5.morfismvmg} 
\pi_{| D |, \; \delta} : V_{| D |, \delta} \longrightarrow 
{\mathcal M}_{g},
\end{equation}where ${\mathcal M}_g$ denotes the moduli space of smooth curves of (geometric) 
genus $g$. Indeed, if 
$F_{\delta}: {\mathcal X}_{\delta} \to V_{| D |, \delta}$ 
denotes the universal family of $\delta$-nodal curves in $S$ 
parametrized by $V_{| D |, \delta}$, 
the fibres of $F_{\delta}$ can 
be simultaneously desingularized, so there exists a diagram of proper 
morphisms
\begin{displaymath}
\begin{array}{cccl}
{\mathcal C}_{\delta} & \stackrel{\Phi_{\delta}}{\longrightarrow} & 
{\mathcal X}_{\delta} &\subset S \times V_{| D |, \delta} \\
 & \searrow^{f_{\delta}} & \downarrow^{F_{\delta}} & \\
 & & V_{| D |, \delta} &    
\end{array} 
\end{displaymath}where 
$\Phi_{\delta}$ is fibrewise the normalization map. In other words $\Phi_{\delta}$ is the blow-up of 
${\mathcal X}_{\delta}$ along its codimension-one singular 
locus and, for each $\delta \geq 1$, 
the morphism $$\pi_{| D |, \delta} : V_{| D |, \delta} 
\to {\mathcal M}_g$$is functorially defined by $f_{\delta}$. 
When $\delta =0$, 
$V_{| D |, 0}$ is the open dense subscheme of smooth curves in 
$|D|$, so $\Phi_0$ is the identity map and we have 
$\pi_{| D |, 0}: V_{| D |, 0} \to {\mathcal M}_{p_a(D)}$. 

\noindent
The problem is to determine, for each morphism $\pi_{| D |, \delta}$, 
the dimension of its image. 

Different from the case of $S= {\Pp}^2$, Severi varieties on surfaces of general type are, in 
general, reducible; for example, Chiantini and Ciliberto (\cite{CC}) showed that even in 
the most natural case of a general surface $S = S_d \subset \Pt$ of degree $d \geq 5$, 
Severi varieties on $S$ of the form $V_{| mH |, \delta}$, $ m \geq d$ and $H$ the 
plane section of $S$, always admit at least one (generically) regular component 
but, sometimes, also some other superabundant components with a dimension bigger than 
the expected one. On the other hand, there are also some results which give upper-bounds 
on $m$ and $\delta$ ensuring that all the components of such a Severi variety are regular 
(see \cite{CS} and \cite{F}). Thus, to precisely approach the problem, we make the following:

\begin{definition}\label{def:nummod}
Let $S$ be a smooth, projective surface of general type and let $D$ be a smooth, irreducible 
curve on $S$. Let $\delta \geq 0$ be such that $V_{| D |, \; \delta} \neq \emptyset$. 
If $V \subseteq V_{| D |, \; \delta}$ is an irreducible component, then the 
{\em number of moduli} of the family of curves parametrized by $V$ is 
$$\nu_{D, \delta} |_V := dim (\pi_{|D|, \delta}(V)).$$
\end{definition}

Since the behaviour of superabundant components is difficult to predict, 
we focus on generically regular components of 
$V_{| D |, \; \delta}$. For this reason, we have to introduce the following 
condition:
\begin{equation}\label{eq:2.restr}
\delta \leq dim(|D|).
\end{equation}Indeed on such a surface, in general, we have 
$dim(|D|) < p_a(D)$ (e.g. if $D$ is a very ample divisor, it directly 
follows from the fact that the characteristic linear system on $D$ is special); therefore 
$V_{| D |, \; \delta}$ cannot have the expected dimension if $\delta >>0$, i.e. if 
$\delta$ is near $p_a(D)$. 

\begin{definition}\label{def:admissible}
The integer $\delta$ will be called {\em admissible} if $\delta$ is as in 
(\ref{eq:2.restr}) and such that $g = p_a(D) - \delta \geq 2$. 
\end{definition}

From Theorem \ref{thm:2.mod.sern} and Remark \ref{rem:0.1} one can eurhystically 
give the following:
\begin{definition}\label{def:5.expmod}
Let $V \subseteq V_{| D |, \; \delta}$ be an irreducible and generically regular component, 
with $\delta$ admissible. Then, the {\em expected number of moduli} of $V$ is 
$$expmod(V) := dim (V).$$ 
\end{definition}

\noindent
Thus, what is expected is that $V$ parametrizes a family having special moduli 
and, moreover, that its number of moduli is the biggest possible; in other words, 
a regular point $[X] \in V \subseteq V_{| D |, \; \delta}$ is expected to be 
birationally isomorphic to finitely many curves in 
$V$. 

By using vector bundle theory on regular surfaces $S$ with effective canonical divisor, one can 
easily determine some examples of regular components of Severi varieties of the form 
$V_{|K_S|, \delta}$ having the expected number of moduli (see \cite{F} for details). On 
the other hand, there are also some 
examples which show that such expected number of moduli is not always 
achieved. Indeed, one can consider particular smooth, projective and regular 
surfaces of general type which belong to a class of 
surfaces that Catanese has recently studied 
(see \cite{Cat1}), called {\em Beauville's surfaces} or 
{\em fake quadrics} (see \cite{Shaf}, page 195). Such a surface is 
of the form $S := (C \times C) /G, $ 
where $C$ is a smooth curve of genus 
$g \geq 2$, $G$ is a finite group acting on each 
factor $C$ and freely acting on the product 
$C \times C$ so that the quotient 
is a smooth surface and the projection $ p : C \times C \to S$ is a topological covering. 
Moreover, if $ |G| = (g - 1)^2$ and if the action of $G$ on $C$ is such that 
$C/ G  \cong {\Pp^1}$, then one determines in $S$ an 
isotrivial rational pencil of smooth curves $\overline{C}$ 
of genus $g$, parametrized by an open dense subset of ${\Pp}^1$. 
From the exact sequence
$$ 0 \to {\Oc}_S \to {\Oc}_S(\overline{C}) 
\to {\Oc}_{\overline{C}}(\overline{C}) \to 0, $$the 
regularity of $S$ and the fact that 
$deg({\Oc}_{\overline{C}}(\overline{C}))=0$, we get 
that $dim (|{\Oc}_S(\overline{C})|) = 1$, so the 
complete linear system coincides with the isotrivial family. Therefore, the morphism 
$$\pi_{|\overline{C}|, 0} : V_{|\overline{C}|, 0} \to {\mathcal M}_{p_a(\overline{C})}$$is 
constant.  

\noindent
{\bf Remark}. The previous example shows that we cannot expect to have, 
"tout court", the expected number of moduli, even in the case of families of 
smooth curves on smooth, projective, minimal and regular surfaces 
of general type. 

From what observed, it is natural to give the following: 

\begin{definition}\label{def:5.moduliproblem}
With the same conditions as in Definition 
\ref{def:5.expmod}, the {\em moduli problem} consists in determining 
for which kind of divisor classes $D \in Div(S)$, 
the number of moduli of generically regular components 
$V \subseteq V_{| D |, \delta}$ coincides with 
the expected one, i.e. when 
$$\nu_{D, \delta} |_V = expmod(V)$$holds.
\end{definition}

Our approach to the moduli problem is analogous to that of Sernesi in \cite{S}, where 
he applied infinitesimal deformation 
theory to families of plane nodal curves. 
This uses the exact sequence (\ref{eq:5.normalmap}). 

\noindent
When, in particular, $V_{| D |, 0}$ 
is considered, if we denote always 
by $X$ the general (smooth) element of $|D|$, then $N = \emptyset$ and 
the Zariski tangent space to $V_{| D |, 0}$ at $[X]$ coincides with 
$H^0(\Oc_S(D))/<X>, $ reflecting the fact that 
$V_{| D |, 0}$ is an open dense subscheme of $|D|$. Moreover, the exact sequence 
(\ref{eq:5.normalmap}) reduces to the standard normal sequence of $X$ in $S$. 
Therefore, if $X$ is a smooth element in $|D|$, we get 
\begin{equation}\label{eq:5.1}
 0 \to H^0({\T}_S |_X) \to H^0({\N}_{X/S}) 
\stackrel{\partial}{\longrightarrow} H^1({\T}_{X}) \to 
\cdots,
\end{equation}where $h^1({\T}_{X}) = 3p_a(X) - 3 = dim( 
{\mathcal M}_{p_a(X)})$, by assumption on $p_a(X) = p_g(X)$.

\noindent
On the other hand, if $[X] \in V_{| D |, \delta}$, $\delta \geq 1$, 
from (\ref{eq:5.normalmap}) we get
\begin{equation}\label{eq:5.3}
 0 \to H^0(\varphi^*({\T}_S)) \to H^0({\N}_{\varphi}) 
\stackrel{\partial}{\longrightarrow} H^1({\T}_{C}) \to 
\cdots, 
\end{equation} where $h^1({\T}_{C}) = 
3g-3 = 3(p_a(X) - \delta -1)= dim( {\mathcal M}_{g})$, with $g \geq 2$ by assumption.

\noindent
Therefore, when $[X] \in V_{| D |, \delta}$, $\delta \geq 0$, is a regular point, the compositions 
$$T_{[X]}(V_{| D |, 0})\hookrightarrow 
H^0({\N}_{X/S}) \stackrel{\partial}{\longrightarrow} H^1({\T}_{X}) $$and 
$$T_{[X]}(V_{| D |, \delta})\hookrightarrow H^0({\N}_{\varphi}) 
\stackrel{\partial}{\longrightarrow} H^1({\T}_{C}), \;\;\; \delta \geq 1,$$can be identified 
with the differentials of the morphisms 
$\pi_{| D |, \delta}$, $\delta \geq 0$, 
at the points $[X]$ and $[C \to X \subset S]$, respectively.

\begin{remark}\label{rem:5.expmod}
\normalfont{If $[X] \in V_{| D | , \delta}$, $\delta \geq 0$, is a regular point and if 
we further assume that $[X]$ is a general point of an irreducible component 
$V$ of $V_{| D | , \delta}$, to give positive answers to the moduli problem of 
Definition \ref{def:5.moduliproblem} we need to show that the differential 
$(\pi_{|D|, \delta})_{*, [X]}$ is injective. From $(\ref{eq:5.1})$ and $(\ref{eq:5.3})$ 
this reduces to finding, for which divisor classes $D$, $h^0(D, {\T}_S|_D) =0$ and 
$h^0(C, \varphi^*({\T}_S))=0$ hold, respectively.}
\end{remark}

\section{The main result}\label{S:3}

From what observed in Remark \ref{rem:5.expmod}, we start by proving the following 
general result.

\begin{theorem}\label{thm:5.Ein} 
Let $S$ be a smooth, projective surface of general type. Let $X \sim D$ be an 
irreducible, $\delta$-nodal curve, $\delta \geq 0$, 
whose set of nodes is denoted by $N$. Then, 

\begin{equation}\label{eq:5.Ein}
h^1(S, \; {\Ii}_{N/S} \otimes \Omega_S^1(D + K_S) ) = 0 \Rightarrow 
h^0(C, \; \varphi^* ({\T}_S))=0.
\end{equation}In particular, when $\delta =0$, 

\begin{equation}\label{eq:5.Ein1}
h^1(S, \; \Omega_S^1(D + K_S)) = 0 \Rightarrow 
h^0(D, \; {\T}_S |_D)=0
\end{equation}
\end{theorem}
\begin{proof}
If $N \neq \emptyset$, denote by $\mu : \tilde{S} \to S$ the blow-up of 
$S$ along $N$, so that one can consider the following diagram of morphisms:
\begin{displaymath}
\begin{array}{clccl}
C &  & \subset & \tilde{S} & \\
\downarrow & \!\!\!\!^{\varphi} & & \downarrow & \!\!\!^{\mu} \\
X & & \subset & S & . 
\end{array}
\end{displaymath}Thus,
$$H^0 (\varphi^*({\T}_S)) = H^0(\mu^*({\T}_S) |_C).$$If we tensor 
the exact sequence defining $C$ in $\tilde{S}$ with $\mu^*({\T}_S)$, 
we get
\begin{equation}\label{eq:5.Einconti}
0 \to \mu^*({\T}_S)(-C) \to \mu^*({\T}_S) \to \mu^*({\T}_S) |_C \to 0.
\end{equation}Observe that
$$H^0 (\mu^*({\T}_S)) \cong H^0({\T}_S)=(0), $$since $H^0({\T}_S)$ 
is isomorphic to the Lie algebra of the Lie group $Aut(S)$, which is finite 
by assumption on $S$ (see \cite{Mat}); thus, the cohomology 
sequence associated to (\ref{eq:5.Einconti}) reduces to 
$$ 0 \to H^0(\mu^*({\T}_S) |_C) \to H^1(\mu^*({\T}_S)(-C)) 
\to \cdots .$$A sufficient condition for 
$h^0(\varphi^*({\T}_S))=0$ is 
therefore $h^1(\mu^*({\T}_S)(-C)) =0$. By Serre duality on $\tilde{S}$, we 
have
\begin{equation}\label{eq:5.Einconti1} 
h^1(\mu^*({\T}_S)(-C)) = h^1((\mu^*({\T}_S))^{\vee} \otimes 
\Oc_{\tilde{S}}(K_{\tilde{S}} + C )).
\end{equation}Since ${\T}_S$ is locally free, 
then $\mu^*({\T}_S)^{\vee} = \mu^*({\T}_S^{\vee}) = 
\mu^*(\Omega^1_S)$, so $(\ref{eq:5.Einconti1})$ becomes
\begin{equation}\label{eq:5.Einconti2} 
h^1(\mu^*({\T}_S)(-C)) = h^1(\mu^*(\Omega_S^1)(K_{\tilde{S}} + C )).
\end{equation}Denote by $B$ the $\mu$-exceptional divisor in $\tilde{S}$ 
such that $ B = \Sigma_{i=1}^{\delta} E_i$. From 
standard computations with blow-ups, we get 
$K_{\tilde{S}} + C = \mu^*(K_S + X) -B.$ Therefore, the right-hand side 
of $(\ref{eq:5.Einconti2})$ becomes $ h^1( \mu^*(\Omega_S^1(K_S 
+ X)) \otimes \Oc_{\tilde{S}}(-B))$. Since we have
$$H^1( \mu^*(\Omega_S^1(K_S + X)) \otimes \Oc_{\tilde{S}}(-B)) \cong 
H^1({\Ii}_{N/S}(X + K_S) \otimes \Omega_S^1),$$from the fact that 
$X \sim D$ on $S$, we get $(\ref{eq:5.Ein})$.

\noindent
For (\ref{eq:5.Ein1}), i.e. $\delta =0$, one can directly use the 
exact sequence
$$0 \to {\T}_S(-D) \to {\T}_S \to {\T}_S|_D \to 0.$$ 
\end{proof}

As an application of Remark \ref{rem:5.expmod} and Theorem \ref{thm:5.Ein}, 
the moduli problem of Definition \ref{def:5.moduliproblem} reduces to 
finding for which 
divisors $D$ on $S \subset {\Pp}^r$ the conditions 
\begin{equation}\label{eq:5.van}
H^1(S, \; \Omega^1_S(K_S + D))= (0)
\end{equation}and 
\begin{equation}\label{eq:5.van1}
H^1(S, \; \Ii_{N/S} \otimes \Omega^1_S(K_S + D))= (0)
\end{equation}hold. The main results of this section 
(Theorems \ref{thm:5.vannodi} and \ref{thm:5.van}) determine 
sufficient conditions on $D$ implying (\ref{eq:5.van}) and (\ref{eq:5.van1}). 

\vskip 4pt

\noindent
{\bf Remark} To prove the basic Lemma \ref{lem:referee} and 
Theorem \ref{thm:5.vannodi}, we shall use some projective-bundle arguments by following the approach 
of \cite{H2}, Sect. II.7. Thus, in the following two results, if $\E$ is a vector bundle on a 
smooth, projective variety $Y$, $\Pp_{Y} ({\E})$ denotes the {\it projective 
space bundle} on $Y$, defined as $Proj(Sym({\E}))$. We have a surjection 
$\pi^*(\E) \to \Oc_{\Pp_{Y} ({\E})}(1)$, where $\Oc_{\Pp_{Y} ({\E})}(1)$ is the 
{\it tautological line bundle} on $\Pp_{Y} ({\E})$ and where $\pi : \Pp_{Y} ({\E}) \to Y$ is 
the natural projection morphism.

\begin{lemma}\label{lem:referee}
Let $S \subset \PR$ be a smooth surface and let $\E$ be a rank 2 vector bundle on $S$. Assume 
that $\E$ is big and nef on $S$ (i.e. the tautological line bundle 
$\Oc_{\Pp_{S} ({\E})}(1)$ is big and nef on $\Pp_{S} ({\E})$). Then 
$$H^i(S, \omega_S \otimes {\E} \otimes det({\E}) \otimes \Oc_S(L)) =(0), $$for $i >0$ and for each 
nef divisor $L$.
\end{lemma}
\begin{proof}By definition, $\Pp_{S} ({\E})$ is a smooth projective variety. From the assumptions 
on $\E$ and $L$ and from the Kawamata-Viehweg vanishing theorem 
(see, for example, \cite{Miy}, page 146), it follows that 
\begin{equation}\label{eq:ref.(*)}
H^i(\Pp_{S} ({\E}), \omega_{\Pp_{S} ({\E})} \otimes 
\Oc_{\Pp_{S} ({\E})}(m) \otimes \pi^*(\Oc_S(L))) = (0), 
\; \; {\rm for} \; i, \; m >0.
\end{equation}

\noindent
Consider the natural projection morphism $\pi : \Pp_{S} ({\E}) \to S$ and recall that 
$$\pi_*(\Oc_{\Pp_{S} ({\E})}(m)) \cong Sym^m({\E}), \; m \geq 1, \; {\rm and} 
\; \pi_*(\Oc_{\Pp_{S} ({\E})}) \cong \Oc_S,$$(see \cite{H2}, 
Prop. II.7.11). From the relative Euler sequence 
$$0 \to \Oc_{\Pp_{S} ({\E})} \to \pi^*({\E}^{\vee}) \otimes \Oc_{\Pp_{S} ({\E})}(1) \to 
{\T}_{\Pp_{S} ({\E})/S} \to 0$$and from the exact sequence
$$0 \to {\T}_{\Pp_{S} ({\E})/S} \to {\T}_{\Pp_{S} ({\E})} \to \pi^*({\T}_S) \to 0, $$we get that 
$$\omega_{\Pp_{S} ({\E})} \cong \Oc_{\Pp_{S} ({\E})}(-2) \otimes 
\pi^*(\omega_S \otimes det({\E})).$$Therefore, if we consider $m=3$ in (\ref{eq:ref.(*)}), we get
\begin{equation}\label{eq:ref.(**)}
H^i(\Pp_{S} ({\E}), \pi^*(\omega_S \otimes det({\E})\otimes \Oc_S(L)) \otimes 
\Oc_{\Pp_{S} ({\E})}(1)) = (0), \; {\rm for} \; i>0.
\end{equation}

\noindent
By projection formula, 
\begin{equation}\label{eq:ref.(***)}
R^i \pi_*(\pi^*(\omega_S \otimes det({\E})\otimes \Oc_S(L)) \otimes 
\Oc_{\Pp_{S} ({\E})}(1)) \cong \omega_S \otimes det({\E}) \otimes \Oc_S(L) 
\otimes R^i\pi_*(\Oc_{\Pp_{S} ({\E})}(1)), 
\end{equation}for each $i >0$. Since the fibres of $\pi$ are isomorphic to 
$\Pp^1$ and since $\Oc_{\Pp_{S} ({\E})}(1)$ is relatively ample, all the higher direct image 
sheaves in (\ref{eq:ref.(***)}) are zero; thus, by Leray spectral sequence and by 
(\ref{eq:ref.(**)}), we get the statement.
\end{proof}

Now, we can prove our main result.

\begin{theorem}\label{thm:5.vannodi}
Let $S \subset {\Pp}^r$ be a smooth surface of general type with hyperplane divisor $H$. 
Suppose that the linear system 
$|D|$ on $S$ has general element which is a smooth, 
irreducible curve. 
Let $X \sim D$ be an irreducible, $\delta$-nodal curve 
of geometric genus $g = p_a(D) - \delta$, where $\delta \geq 0$ admissible (as in 
Definition \ref{def:admissible}). 
Assume that:
\begin{itemize}
\item[(i)] $\Omega^1_S(K_S)$ is globally generated;
\item[(ii)] $D \sim K_S + 6H + L$, where $L$ is a nef divisor;
\item[(iii)] the Severi variety $V_{| D |, \delta}$ is regular 
at $[X]$ (in the sense of Definition \ref{def:2.reg}). 
\end{itemize}Then, the morphism 
$$\pi_{| D |, \delta} : V_{| D |, \delta} \to {\mathcal M}_g $$has 
injective differential at $[X]$. In particular, 
$\pi_{| D |, \delta} $ 
has finite fibres on each generically regular component of $V_{| D |, \delta}$, 
so each such component parametrizes a family 
having the expected number of moduli.
\end{theorem}
\begin{proof}
First of all, we want to show that hypothesis (ii) implies 
(\ref{eq:5.van}). To prove this, we will use Lemma \ref{lem:referee}. 
Therefore, the first step of our analysis is 
to apply such vanishing result to the vector bundle 
$${\E} = \Omega_S^1(aH),$$where 
$a$ is a positive integer. The problem reduces to finding which "twists" of 
$\Omega_S^1$ are big and nef 
on $S \subset {\Pp}^r$. In the sequel we shall write for short 
$\Omega_S^1(a)$ instead of $\Omega^1_S(aH)$. From the exact sequence
$$0 \to {\it Con}_{S/{\Pp}^r}(a) \to \Omega^1_{{\Pp}^r}(a)|_S 
\to \Omega^1_S(a) \to 0,$$it is useful 
compute for which positive integers $a$ the vector bundle 
$\Omega^1_{{\Pp}^r}(a)$ 
is ample or globally generated (see \cite{H1}). From the Euler sequence of 
$\PR$ one deduces the exact sequence
$$ 0 \to \Omega^2_{{\Pp}^r} \to {\Oc}_{{\Pp}^r}^{\oplus \frac{r(r+1)}{2}}(-2) 
\to \Omega^1_{{\Pp}^r}\to 0$$(see \cite{OSS}, page 6, 
and \cite{SS}, page 73); therefore, one trivially has
$$ 0 \to \Omega^2_{{\Pp}^r}(2)\to {\Oc}_{{\Pp}^r}^{\oplus \frac{r(r+1)}{2}} 
\to \Omega^1_{{\Pp}^r}(2)\to 0,$$i.e. 
$\Omega^1_{{\Pp}^r}(2)$ and, so, $\Omega^1_S(2)$ are globally generated whereas 
$\Omega^1_{{\Pp}^r}(a)$ and $\Omega^1_S(a)$ are ample, for $a \geq 3$. 

Recall now that $\Pp_{\PR} (\Omega^1_{{\Pp}^r}(1))$ is the universal line over the Grassmannian 
${\Gr}(1,r)$ of lines in $\PR$ (see, for example, \cite{J}, app. B and C, or \cite{K}, page 
369). By standard properties of projective bundles, $\Pp_{\PR} (\Omega^1_{{\Pp}^r}(1)) \cong 
\Pp_{\PR} (\Omega^1_{{\Pp}^r}(2))$, thus we have 
$${\Ff} := \Pp_{\PR} (\Omega^1_{{\Pp}^r}(2)) \subset {\Gr}(1,r) \times \PR$$with the natural 
projection $p_i$ on the $i$-th factor, $1 \leq i \leq 2$. If $\gamma$ denotes 
the Pl\"ucker embedding of ${\Gr}(1,r)$ in ${\Pp}^{\oplus \frac{r(r+1)}{2}-1} $, one determines the map 
$$f := \gamma \circ p_1 : {\Ff} \to {\Pp}^{\oplus \frac{r(r+1)}{2}-1}.$$On the other hand, 
we can consider the complete tautological linear system $|\Oc_{\Ff}(1)|$, which is free since 
$\Omega^1_{\PR}(2)$ is globally generated. From the Leray spectral sequence, the Euler sequence 
and the Bott formula (see \cite{OSS}, page 8), we get that
$$H^0({\Ff}, \Oc_{\Ff}(1)) \cong H^0(\PR , \Omega^1_{\PR}(2)) \cong \bigwedge^2V,$$where here 
$\PR = \Pp(V) = Proj(Sym(V))$. Therefore, the complete linear system $|\Oc_{\Ff}(1)|$ 
defines a morphism 
$$\Phi : {\Ff} \to \Pp(\bigwedge^2 V^*) \cong {\Pp}^{\oplus \frac{r(r+1)}{2}-1}.$$One easily sees that 
$\Phi$ and $f$ coincide, so the global sections of 
$\Oc_{\Ff}(1)$ contract the $p_1$-fibres of ${\Gr}(1,r)$ in ${\Ff}$, which 
are lines in $\PR$. 

\noindent
From the fact that 
$\Pp_{S} (\Omega^1_S(2)) \subset \Pp_{S} (\Omega^1_{\PR}|_S(2))$, 
the restriction of $\Phi$ to $\Pp_{S} (\Omega^1_S(2))$ is generically finite since 
$S$, being of general type, is not filled by lines. Thus the rank 2 vector bundle $\Omega^1_S(2)$ 
is globally generated and big and nef. By Lemma \ref{lem:referee}, 
$H^1(S , \; \omega_S \otimes \Omega^1_S(2) 
\otimes det (\Omega^1_S(2)) \otimes \Oc_S(L)) =(0)$, for each nef divisor $L$. Since 
$det(\Omega^1_S(2)) = \Oc_S(K_S + 4H)$, we have that $H^1(S, \Omega^1_S(2K_S + 6H +L)) = (0),$ for 
each nef divisor $L$.

Therefore, if $D \sim K_S + 6H + L$, with $L$ nef, then 
$$(*) \;\; \; H^1(S, \; \Omega^1_S(K_S + D))=(0).\footnote{Observe that if 
one directly applies Griffiths vanishing results, i.e. Theorem (5.52), Theorem (5.64) and 
Corollary (5.65) in \cite{SS}, to the vector bundle $\Omega^1_S(a)$, $a \geq 2$, one determines 
stronger conditions on $D$. Precisely, $L$ must be ample instead of nef. Therefore, the 
approach above determines more general conditions on $D$.}$$

The vanishing result $(*)$ is a fundamental 
tool for the following second part of the proof. On $S$ we can consider the exact sequence
\begin{equation}\label{eq:5.conto}
0 \to \Ii_{N/S}(D) \to \Oc_S(D) \to \Oc_N(D) \to 0
\end{equation}which determines 
the restriction map $\rho_D$:
$$0 \to H^0(\Ii_{N/S}(D)) \to H^0(\Oc_S(D)) \stackrel{\rho_D}{\rightarrow} 
H^0(\Oc_N(D)) \to H^1(\Ii_{N/S}(D)) \to \cdots .$$By hypothesis (iii), 
$\rho_D$ is surjective. Next, by tensoring the exact sequence 
$(\ref{eq:5.conto})$ with 
$\Omega_S^1(K_S)$, we get
$$0 \to H^0(\Ii_{N/S}(D) \otimes \Omega_S^1(K_S)) \to 
H^0(\Omega_S^1(K_S+D)) 
\stackrel{\rho_{\Omega_S^1(K_S+D)}}{\longrightarrow}$$ 
$$ \stackrel{\rho_{\Omega_S^1(K_S+D)}}{\longrightarrow} 
H^0(\Oc_N(\Omega_S^1(K_S+D))) \cong {\CC}^{2 \delta} \to 
H^1(\Ii_{N/S}(D) \otimes \Omega_S^1(K_S)) \to $$ 
$$\to H^1(\Omega_S^1(K_S+D))\to 0.$$Thus, 
the map $\rho_{\Omega_S^1(K_S+D)}$ 
is surjective if and only if $H^1(\Ii_{N/S}(D) \otimes \Omega_S^1(K_S)) 
\cong H^1(\Omega_S^1(K_S+D))$. From the 
first part of this proof, hypothesis (ii) implies that 
$H^1(\Omega_S^1(K_S+D)) = (0)$, so we have
$$h^1(\Ii_{N/S}(D) \otimes \Omega_S^1(K_S))=0 \Leftrightarrow \; 
\rho_{\Omega_S^1(K_S+D)} \; {\rm surjective}.$$By ($\ref{eq:5.Ein}$) of 
Theorem \ref{thm:5.Ein}, the surjectivity of $\rho_{\Omega_S^1(K_S+D)}$ 
implies therefore that $h^0(\varphi^*({\T}_S))=0 $ and so the statement.

The last step is to determine if, with the given hypotheses, 
the map $\rho_{\Omega_S^1(K_S+D)}$ is surjective. Consider the map
\begin{equation}\label{eq:5.conto1}
H^0(\Omega_S^1(K_S+D)) \stackrel{\rho_{\Omega_S^1(K_S+D)}}{\longrightarrow} 
H^0(\Oc_N(\Omega_S^1(K_S+D)) \cong {\CC}^{2 \delta} \cong 
\bigoplus_{i=1}^{\delta} \CC^2_{(i)}.
\end{equation}By hypothesis (i), for each $p \in S$, the sheaf morphism 
$$H^0(\Omega_S^1(K_S)) \otimes {\Oc}_{S,p} \to \Omega_S^1(K_S)|_p 
\cong {\Oc}_{S,p}^{\oplus 2}$$is surjective; thus, for each 
$p \in S$ there exist two global sections $s_1^p, \; s_2^p \in 
H^0 (\Omega_S^1(K_S))$ which generate the stalk $\Omega^1_S(K_S) |_p$ 
as an ${\Oc}_S$-module, i.e. 
$$s_1^p(p) = (1,0) \;
{\rm and} \; s_2^p(p) = (0,1) \in {\Oc}_{S,p}^{\oplus 2}\; .$$If 
$N = \{p_1, \; p_2, \; \ldots, \;p_{\delta} \}$ is the set of nodes 
of $X$, then $H^0({\Oc}_N(D)) \cong \CC^{\delta} \cong \CC_{(1)} 
\oplus \CC_{(2)} \oplus \cdots \oplus \CC_{(\delta)}$. The surjectivity 
of $\rho_D$ implies there exist global sections $\sigma_i \in 
H^0({\Oc}_S(D))$ 
such that
$$\sigma_i(p_j) = (0, \; 0, \; \ldots, \;0), \; if \; 1 \leq i \neq j \leq 
\delta,$$ 
\vskip -0,6 cm
$$\sigma_i(p_i) = (0, \; \ldots, \; 0, \; 1 , \; 0, \; \ldots, \; 0), 
\; 1 \in \CC_{(i)}, \; 1 \leq i \leq \delta .$$Therefore, 
$s_1^{p_i} \otimes \sigma_i, \; s_2^{p_i} \otimes \sigma_i 
\in H^0(\Omega^1_S(D+K_S))$ and 
$$s_1^{p_i} \otimes \sigma_i(p_j) = s_2^{p_i} \otimes \sigma_i(p_j) 
= (0, \; \ldots, \; 0) \in \CC^2_{(1)}\oplus \cdots \oplus 
\CC^2_{(\delta)} \cong \CC^{2 \delta}, \; 1 \leq i \neq j \leq \delta,$$
\vskip -0,6 cm
$$s_1^{p_i} \otimes \sigma_i(p_i)= ((0 ,0) , \; \ldots, \; (1,0) , 
\; \ldots, \; (0,0)) = (0, \; \ldots, \; 1, \; 0, \; \ldots, \; 
0) \in \CC^{2 \delta},$$where $(1,0) \in \CC_{(i)}^2$ and    
$$s_2^{p_i} \otimes \sigma_i(p_i)= ((0 ,0) , \; \ldots, \; (0,1) , 
\; \ldots, \; (0,0)) = (0, \; \ldots, \; 0, \; 1, \; \ldots, \; 
0) \in \CC^{2 \delta},$$ where $(0,\; 1) \in \CC_{(i)}^2,$ for 
$ 1 \leq i \leq \delta$. This means that the map $(\ref{eq:5.conto1})$ 
is surjective. Moreover, since the condition for a point 
$[X] \in V_{| D |, \delta}$ to be regular is an open condition 
in the family, it follows that the component of 
$V_{| D |, \delta}$ containing $[X]$ has 
the expected number of moduli. 
\end{proof}

From the first part of the proof of Theorem \ref{thm:5.vannodi} we observe 
that in the case of familes of smooth curves one can eliminate 
hypotheses (i) and (iii). Indeed, we have the following result.

\begin{theorem}\label{thm:5.van}
Let $S\subset \PR$ be a smooth surface of general type and let $D$ be an 
effective divisor on $S$. Denote by $H$ the hyperplane 
section of $S$. Assume that $$D \sim K_S + 6 H + L,$$ 
where $L$ is a nef divisor on $S$. Then, $H^1(S, \;\Omega^1_S(K_S + D)) =(0)$. 

\noindent
If, moreover, $|D|$ contains 
smooth, irreducible elements, the family of 
smooth curves $V_{| D |, 0}$ has the expected number of moduli.
\end{theorem}
\begin{proof}
For the first part of the statement, one can repeat the procedure 
at the beginning of the proof of Theorem \ref{thm:5.vannodi}. From 
$(\ref{eq:5.Ein1})$ we get the second part of the statement.
\end{proof}

Let $S = S_d \subset \Pt$ be a smooth surface of 
degree $d$; in view of the fact that $K_S \sim (d-4) H$, as a 
corollary of Theorems \ref{thm:5.vannodi} and \ref{thm:5.van} we get:

\begin{coroll}\label{cor:5.vannodi}
If $S_d \subset {\Pp}^3$ is a smooth surface of degree $d \geq 6$, 
the generically regular components of $V_{| mH |, \delta}$ 
have the expected number of moduli, when $m \geq d+2$ and $\delta \geq 1$ admissible. 
The same conclusion holds for the family of smooth curves $V_{| mH |, 0}$, when 
$d \geq 5$ and $m \geq d+2$. 
\end{coroll}

\begin{remark}\label{rem:rem3.0}
\normalfont{
More generally, if $S \subset \PR$ is of general type 
with $K_S \sim \alpha H$, then 
we have positive answers to the moduli problem for all generically 
regular components of $V_{| mH |, \delta}$, with $m \geq \alpha + 6$, 
when $\alpha \geq 2$ and $\delta \geq 1$ admissible, and with $m \geq \alpha + 6$, when 
$\delta =0$ and $\alpha \geq 1$.   
}
\end{remark}

\begin{remark}
\normalfont{The conditions $\Omega_S^1(2)$ globally generated and big and nef on $S$ 
play a crucial role in the proof of Theorem \ref{thm:5.vannodi}. 
Thus with this approach this result is, in a certain sense, sharp. 
For example, if we focus on regular surfaces, $\Omega_S^1$ cannot be globally generated, since 
$\overline{H^{1,0}}(S) = H^{0,1}(S)$ and $H^{1,0}(S) \cong H^0(S, \Omega_S^1)$ whereas 
$H^{0,1}(S) \cong H^1(S, {\Oc}_S)=(0)$. If $S$ is also a non-degenerate complete intersection 
in $\PR$, then $\Omega^1_S(1)$ cannot be globally generated. 
Furthermore, we have some results of Schneider (see \cite{Schn}) 
which state that, even in the most natural case of smooth surfaces $S_d \subset {\Pp}^3$ 
of degree $d \geq 5$, $\Omega_{S_d}^1$ and $\Omega_{S_d}^1(1)$ are not ample.}
\end{remark}

\section{A fundamental proposition}\label{S:4}
The aim of this and the following two sections is to find other results 
giving positive answers to the moduli problem, posed in Definition 
\ref{def:5.moduliproblem}, for some other classes of divisors on $S$ which are not covered 
by Theorems \ref{thm:5.vannodi} and \ref{thm:5.van}.

From now on we shall focus on the case 
of regular surfaces; therefore $S$ will always denote a smooth, regular surface of general type, 
unless otherwise specified. In such a case, we are able, in particular, to complete 
Remark \ref{rem:rem3.0} by also including 
divisors $D \sim m H$ with $ 1 \leq m \leq \alpha + 5$ and with some further conditions on 
$\delta$.

The first step of our analysis is based on a key proposition concerning 
first-order deformations of the normalization morphism $\varphi: C \to X \subset S$. Then 
we conclude, in some cases, by using a detailed analysis of the Brill-Noether map of the 
line bundle ${\Oc}_C(\varphi^*(H))$, in some other cases, by using uniqueness results of certain 
linear systems on $C$. 

The core of this section is to prove such a fundamental proposition. Before doing this, we need 
to remind some general facts. 

Let $S \subset \PR$ be a smooth, non-degenerate surface (not necessarily 
regular and of general type). As in (\ref{eq:1bis}), 
the normalization morphism $\varphi$ is a map from $C$ to $S$ such that 
$Im(\varphi) = X \subset S$. If $ i : S \hookrightarrow \PR$ is the natural embedding, 
we have the following diagram of morphisms:
\begin{displaymath}
\begin{array}{ccc}
C & &  \\
\downarrow^{\varphi} &  \searrow^{\psi} & \\
S & \stackrel{i}{\hookrightarrow} & \PR
\end{array}
\end{displaymath}where $\psi = i \circ \varphi : C \to \PR$. By pulling back to $C$ the 
normal sequence of $S$ in $\PR$, we get the exact sequence of vector bundles 
on $C$
$$ 0 \to \varphi^*({\T}_S) \to  \psi^*({\T}_{\PR}) \to \varphi^*({\N}_{S/\PR}) \to 0.$$Thus, 
\begin{equation}\label{eq:3.*}
0 \to H^0( \varphi^*({\T}_S)) \to H^0(\psi^*({\T}_{\PR})) 
\to H^0(\varphi^*({\N}_{S/\PR})) \to \cdots
\end{equation}holds, where $H^0( \varphi^*({\T}_S))$ parametrizes 
first-order deformations of the map $\varphi : C \to S$, 
with $C$ and $S$ both fixed, as well as $H^0(\psi^*({\T}_{\PR}))$ parametrizes 
first-order deformations of the map $\psi : C \to \PR$, 
with $C$ and $\PR$ both fixed (see \cite{HM}). 

We also recall the following useful definition (see \cite{CS}).

\begin{definition}\label{def:1}
Let $X$ be any reduced, irreducible curve in $\mathbb{P}^r$. $X$ is said to be 
\emph{geometrically linearly normal} (\emph{g.l.n.} for short) if 
the normalization map $\varphi: C \to X \subset \PR$ cannot be factored 
into a non-degenerate morphism $C \to \mathbb{P}^N$, with $N>r$,
followed by a projection.
\end{definition}

\noindent
In other words, if $H$ is the hyperplane section of $X$, 
$|{\Oc}_C(\varphi^*(H))|$ must be complete.

We are now able to give the following:

\begin{proposition}\label{prop:chiave}
Let $S \subset \PR$ be a smooth, regular, non-degenerate and linearly normal surface of 
general type. Let $[X] \in V_{| D |, \delta}$ be a regular point of the Severi variety 
$V_{| D |, \delta}$ on $S$. 

\noindent
(i) Assume that $X$ is non-degenerate in $\PR$ and geometrically 
linearly normal. If $h^0(C, \psi^*({\T}_{\PR})) = (r+1)^2 -1$, then all 
first-order deformations of the map $\psi: C \to \PR$, with $C$ fixed, are induced 
by first-order projectivities (i.e. by elements of $H^0({\T}_{\PR}))$. 
Moreover, $h^0(\varphi^*(\T_S)) = 0$.

\noindent
(ii) Assume that $ D \sim H$ on $S$ and that $X \subset H \cong {\Pp}^{r-1}$ is 
non-degenerate and g.l.n. as a curve in ${\Pp}^{r-1}$. Suppose also that $S$ is such that 
$h^1({\Oc}_S(H)) = 0$ and $|K_S| \neq \emptyset$. If $h^0(\psi^*(\T_{\PR})) =r^2 + r -1$, 
then all first-order deformations of the map $\psi: C \to \PR$, with $C$ 
fixed, are induced by first-order projectivities not fixing pointwise the hyperplane 
$H \subset \PR$. Moreover, $h^0(\varphi^*(\T_S)) = 0$.
\end{proposition}
\begin{proof} (i) The first part of the statement is a straightforward computation. We shall 
briefly recall the fundamental steps of its proof. If $\mu : \tilde{S} \to S$ is 
the blow-up of $S$ along $N = Sing(X)$, by the 
hypotheses on $S$ and by the pull-back to $\tilde{S}$ of the Euler sequence, we get
$$H^0({\T}_{\PR}) \cong H^0({\T}_{\PR}|_S) \cong H^0(\mu^*({\T}_{\PR})).$$Since $X$ is 
g.l.n. and non-degenerate in $\PR$, by Serre duality and by the 
pull-back on $C$ of the Euler sequence we have
$$(*) \;\;\; 0 \to K \to H^0(\psi^*({\T}_{\PR})) \to (coker(\mu_{0,C}))^{\vee} \to 0,$$where 
$K=  (H^0 ({\Oc}_C(\varphi^*(H)))^{\vee} \otimes H^0({\Oc}_C(\varphi^*(H))))/ 
H^0({\Oc}_C) \cong H^0({\T}_{\PR})$ and 
where$$\mu_{0,C} : H^0({\Oc}_C(\varphi^*(H)) \otimes 
H^0 (\omega_C(-\varphi^*(H))) \to H^0(\omega_C)$$is the Brill-Noether map of 
${\Oc}_C(\varphi^*(H))$. Since 
$h^0({\T}_{\PR}) = dim (PGL(r+1, \CC)) = (r+1)^2 -1$, from $(*)$ it follows that 
$h^0(\psi^*({\T}_{\PR})) = (r+1)^2 - 1$ iff $dim (coker(\mu_{0,C})) = 0$. In this case, 
by standard 
Brill-Noether theory (see \cite{ACGH}, Proposition 4.1, 
page 187), there is no first-order 
deformation of $\psi: C \to \PR$, with $C$ fixed, induced by 
first-order deformations of the linear system 
$|{\Oc}_C(\varphi^*(H))|$; so all such deformations are induced by elements of 
$H^0({\T}_{\PR})$. 

To get the second part of statement (i) observe that, by the regularity of $S$ and 
by (\ref{eq:2.Y2}), 
$H^0({\N}_{\varphi}) \cong T_{[X]}(V_{|D|, \delta})$. Assume, 
by contradiction, that 
$h^0(\varphi^*({\T}_S))\neq 0$ and let $ v \in H^0(\varphi^*({\T}_S))$ be a non-zero 
vector. Such $v$ corresponds to a tangential direction 
$v \in T_{[X]}(V_{|D|, \delta})$ 
since, by (\ref{eq:5.3}), we have 
$H^0(\varphi^*({\T}_S) ) \subseteq H^0({\N}_{\varphi})$. 

By the regularity assumption of $[X] \in V_{|D|, \delta}$, all directions in 
$T_{[X]}(V_{|D|, \delta})$ are unobstructed. This means there exist a one-dimensional base 
scheme $\Delta$, smooth at the central point $o \in \Delta$, and a family 
${\mathcal X} \to \Delta$ such that 
$${\mathcal X} = \{ X_t \}_{t \in \Delta} \subset S \times \Delta$$where 
$$[X_t] \in V_{|D|, \delta}, \; \forall \; t \in \Delta, \; 
[X_o] = [X] , \; {\rm and} \; T_{[o]} (\Delta) = 
< v >. $$Since $ < v > \subset H^0({\varphi}^*({\T}_S))$, the 
family ${\mathcal X} \to \Delta$ corresponds to a family of maps 
$\Phi : C \times \Delta \to S \times \Delta$, for which 
$$\Phi = \{ \varphi_t \}_{t \in \Delta}, \; \varphi_o = \varphi, \; 
\varphi_t = \Phi|_t : C \times \{ t \} \to S \times \{ t \}, \; 
\varphi_t(C) = X_t \subset S.$$By composing $\Phi$ with the map 
$i \times id_{\Delta}$, where $i : S \hookrightarrow \PR$, we get a family of maps 
$\Psi : C \times \Delta \to {\PR} \times \Delta$ for which  
$$\Psi = \{ \psi_t \}_{t \in \Delta}, \; \psi_o = \psi, \; 
\psi_t = \Psi|_t : C \times \{ t \} \to {\PR} \times \{ t \}, \; 
\psi_t(C) = X_t \subset S \subset \PR. $$From (\ref{eq:3.*}), we know that 
$H^0 (\varphi^*({\T}_S)) \subseteq H^0 (\psi^*({\T}_{\PR}))$ and, from the above 
computations, we have $H^0( \psi^*({\T}_{\PR})) \cong H^0 ({\T}_{\PR})$. Therefore, 
the element $v \in  H^0( \psi^*({\T}_{\PR}))$ is induced by first-order 
projectivities, so the family $\Psi \to \Delta$ is determined by a family 
$\Omega \to \Delta$, where
$$\Omega \subset PGL(r+1, \CC), \; \Omega: X \times \Delta \to {\PR} \times \Delta, \; 
\Omega = \{ \omega_t \}_{t \in \Delta}$$such that 
$\psi_t = \omega_t \circ \psi$ and 
$[\omega_t(\psi(C))] = [X_t] \in V_{|D|, \delta}$, for each $t \in \Delta$ whereas 
$[\omega_o(\psi(C))] = [X]$. 

Since $S$ is of general type, then $\Omega \subset PGL(r+1, \CC) \setminus Aut(S)$. Therefore, 
if$$ X_t = \omega_t(X) \subset S, \; \forall \; t \in \Delta,$$then
$$X \subset \omega_t^{-1}(S) = S_t, \; \forall \; t \in \Delta,$$where $S_t \subset \PR$ 
is a smooth surface projectively equivalent to $S$, for each $t \in \Delta$, and 
$S_o = S$. We therefore obtain a family of maps $\Lambda: S \times \Delta  
\to {\PR} \times \Delta$ such that 
$\Lambda |_t = \omega^{-1}_t$, for each $t \in \Delta$. By composing such family of maps 
with $\mu \times id_{\Delta}$, $\mu : {\tilde S} \to S$, we thus get a family of maps 
$$\Theta : {\tilde S} \times \Delta \to {\PR} \times \Delta$$where 
$$\Theta |_t = \omega_t^{-1} \circ \mu : {\tilde S} \times \{t \} \to {\PR} \times \{t \}, 
\; \Theta |_t ({\tilde S}) = \omega_t^{-1}(\mu ({\tilde S}))= \omega_t^{-1}(S) = S_t \subset 
\PR.$$Since $\Theta |_o = id_{\PR} \circ \mu$ and $T_{o} (\Delta) = < v>$, 
the element $v \in H^0(\varphi^*({\T}_S)) \subset H^0 (\psi^*({\T}_{\PR}))$ is also an 
element of $H^0(\mu^*({\T}_{\PR}) \otimes {\Oc}_{\tilde{S}}(-C))$. 

This leads to a contradiction; indeed, by tensoring the exact sequence defining $C$ in 
${\tilde S}$ with $\mu^*({\T}_{\PR})$, we get
$$0 \to \mu^*({\T}_{\PR}) \otimes {\Oc}_{\tilde{S}}(-C)  \to  \mu^*({\T}_{\PR}) \to 
 \mu^*({\T}_{\PR})|_C \cong \psi^*({\T}_{\PR}) \to 0.$$From the above computations, 
we know that $H^0(\mu^*({\T}_{\PR})) \cong H^0(\psi^*({\T}_{\PR}))$, which 
implies $h^0(\mu^*({\T}_{\PR}) \otimes {\Oc}_{\tilde{S}}(-C)) = 0$. 

\vskip 0.2cm

\noindent
(ii) In this case $X \sim H$ on $S$ and $X \subset H \cong {\Pp}^{r-1}$ is non-degenerate 
in ${\Pp}^{r-1}$, then
$$\psi^*({\T}_{\PR}) \cong \psi^*({\T}_{{\Pp}^{r-1}}) \oplus 
{\Oc}_{C}(\psi^*(H))$$(with abuse of notation, we denote always by $\psi$ the map 
$\psi: C \to X \subset H \cong {\Pp}^{r-1}$). From the hypotheses on $X$, we get 
$$h^0(\psi^*({\T}_{\PR})) = h^0(\psi^*({\T}_{{\Pp}^{r-1}})) + r.$$By using the same computations 
of (i), we get$$ 0 \to H^0({\Oc}_C) \to  H^0 ({\Oc}_C(\psi^*(H))^{\vee} \otimes 
H^0({\Oc}_C(\psi^*(H)) \to  H^0(\psi^*({\T}_{{\Pp}^{r-1}})) \to (coker(\mu_{0,C}))^{\vee} 
\to 0,$$where
$$(coker(\mu_{0,C}))^{\vee}  \cong T_{[|{\Oc}_C(\psi^*(H))|]}(G^{r-1}_{deg(X)}(C)).$$Thus, 
as in (i), $h^0(\psi^*({\T}_{{\Pp}^{r-1}})) = r^2 -1$ if and only if 
$dim(coker(\mu_{0,C}))=0$. 

Note that
\begin{equation}\label{eq:prop4.1*}
\frac{H^0({\tilde S},  \mu^*({\T}_{\PR}))}{H^0({\tilde S}, 
\mu^*({\T}_{\PR}) \otimes {\Oc}_{\tilde{S}}(-C)))} \stackrel{\beta}{\hookrightarrow} 
H^0(C, \psi^*({\T}_{\PR})).
\end{equation}

From the pull-back of the Euler sequence and from the hypotheses on $S$, we get
\begin{equation}\label{eq:prop4.1***}
0 \to {\Oc}_{\tilde{S}}(-C)  \to H^0({\Oc}_{\tilde{S}}(\mu^*(H)))^{\vee}\otimes 
{\Oc}_{\tilde{S}}(\mu^*(H) - C)  \to \mu^*({\T}_{\PR}) \otimes {\Oc}_{\tilde{S}}(-C) 
\to 0 .
\end{equation}Observe that $h^0({\Oc}_{\tilde{S}}(-C)) = h^1({\Oc}_{\tilde{S}}(-C)) = 0$: indeed, 
the first vanishing trivially holds whereas, by Leray's isomorphism and by 
Serre duality, we have $h^1(\Oc_{\tilde{S}}(-C))= h^1(\Ii_{N/S}(K_S + H))$; from the regularity 
of $[X] \in V_{|H|, \delta}$, Remark \ref{rem:1.4} and the hypothesis $h^1({\Oc}_S(H)) = 0$, 
we get $ h^1(\Ii_{N/S}( H)) = 0$. Since $K_S$ is effective by assumption, 
$N$ also imposes independent conditions to $|K_S + H|$. By standard Mumford's vanishing 
theorem, we have $h^1({\Oc}_S(K_S + H)) = 0$, so $h^1(\Ii_{N/S}(K_S + H)) = 0$. 

\noindent
We therefore obtain
\begin{displaymath}
\begin{array}{cl}
H^0(\mu^*({\T}_{\PR}) \otimes {\Oc}_{\tilde{S}}(-C)) & \cong 
H^0 ({\Oc}_{\tilde{S}}(\mu^*(H)))^{\vee}\otimes 
H^0({\Oc}_{\tilde{S}}(\mu^*(H) - C)) \\
   & = H^0 ({\Oc}_{\tilde{S}}(\mu^*(H)))^{\vee}\otimes H^0({\Oc}_{\tilde{S}}(2B)),
\end{array}
\end{displaymath}where $B = \Sigma_{i=1}^{\delta} E_i$ is the 
$\mu$-exceptional divisor. Since $2B$ is a fixed divisor, 
$h^0(\mu^*({\T}_{\PR}) \otimes {\Oc}_{\tilde{S}}(-C)) = 
h^0 ({\Oc}_{\tilde{S}}(\mu^*(H)) = r+1$. Moreover, since 
$H^0(\mu^*({\T}_{\PR}) \otimes {\Oc}_{\tilde{S}}(-C)) \subset 
H^0(\mu^*({\T}_{\PR})) \cong H^0({\T}_{\PR})$, the elements of such a vector space correspond 
to first-order projectivities fixing pointwise the hyperplane $H \subset \PR$. 
Turning back to (\ref{eq:prop4.1*}), $h^0(\psi^*({\T}_{\PR})) = r^2 + r -1$ if and only if 
$\beta$ is an isomorphism. In such a case, all first-order deformations of 
$\psi: C \to \PR$, with $C$ fixed, are induced up to first-order 
by projectivities not fixing pointwise the curve $X \subset H$.

For the second statement in (ii), one can follow the same procedure in (i). By supposing there 
exists a non-zero vector $v \in H^0(\varphi^*({\T}_S))$, one determines a 
family $\Omega \to \Delta$, where $\Omega \subset PGL(r+1, \CC) \setminus Aut(S)$, such that
$$\Omega = \{ \omega_t \}_{t \in \Delta}, \; \omega_t(X) = X_t \subset S, \; {\rm and} 
\; T_o(\Delta) = <v>.$$As before, one obtains 
$v \in H^0(\mu^*({\T}_{\PR}) \otimes \Oc_{\tilde{S}}(-C))$, so 
the family $\Omega$ is contained in the sugroup 
$\Gamma < PGL(r+1, \CC)$, whose elements pointwise fix the curve $X$. Therefore, 
we have $\omega_t(X) = X$, for each $t \in \Delta$, 
contradicting the existence of the non-trivial, one-dimensional family ${\mathcal X} = 
\{ X_t \}_{t \in \Delta}$.
\end{proof}  

From Remark \ref{rem:5.expmod}, in the sequel we will be concerned in 
finding conditions which imply the hypotheses of Proposition \ref{prop:chiave}. 
These will give further affirmative answers to the moduli problem posed in Definition 
\ref{def:5.moduliproblem} for Severi varieties on smooth, regular and non-degenerate surfaces 
$S\subset \PR$ of general type.

\section{Number of moduli for families of non-degenerate, nodal curves on 
linearly normal surfaces of general type.}\label{S:5}

As remarked at the beginning of Section \ref{S:4}, we want to find some other conditions 
establishing positive answers to the moduli problem for those Severi varieties which do not 
satisfy the hypotheses of Theorems \ref{thm:5.vannodi} and \ref{thm:5.van}. 

Here we shall focus on the case of $S \subset \PR$ a smooth surface of 
general type which is regular, non-degenerate, linearly normal and such that 
$h^1({\Oc}_S(H)) = 0$, $H$ the hyperplane section of $S$. 
Observe that, in this case, one can obviously apply the 
results in Section \ref{S:3}, since they are more generally valid.

The results we obtain here apply, for example, to some cases which are not covered 
by Corollary \ref{cor:5.vannodi} and Remark \ref{rem:rem3.0} even though their statement 
gives some restrictions to the 
admissible number of nodes $\delta$ with respect to (\ref{eq:2.restr}).  

In this section, we consider $[X] \in V_{| D |, \delta}$ on $S$ 
such that $X$ is non-degenerate in $\PR$. 
From Proposition \ref{prop:chiave} (i), we want to find conditions on $D$ 
in order that $h^0(\psi^*(\T_{\PR})) = (r+1)^2 -1 = dim(PGL(r+1, \CC))$. 
To this aim, put 
\begin{equation}\label{eq:19bis}
{\Oc}_C(\psi^*(H)) = {\Oc}_C(\tilde{H}),
\end{equation}then 
$X$ is geometrically linearly normal (see Definition \ref{def:1}) 
if and only if $|{\Oc}_C(\tilde{H})|$ 
is complete of dimension $r$. In such a case, we consider the 
{\em Brill-Noether map} of the line bundle ${\Oc}_C(\tilde{H})$, i.e. 
\begin{equation}\label{eq:BNmap}
\mu_{0,C} : H^0({\Oc}_C(\tilde{H})) \otimes 
H^0 (\omega_C(-\tilde{H})) \to H^0(\omega_C).
\end{equation}

\begin{remark}\label{rem:pareschi}
\normalfont{
Similarly to Definition 1.1.2 in \cite{Par}, 
if $X$ is g.l.n. and if the map $\mu_{0,C} $ is surjective, then 
$| {\Oc}_C(\tilde{H}) |$ is called an {\em isolated linear system} on $C$. 
The surjectivity of $\mu_{0,C} $ implies the 
injectivity of the dual map $\mu_{0,C}^{\vee}$ 
so the Euler exact sequence on $C$,
\begin{equation}\label{eq:Euler}
0 \to {\Oc}_C \to 
H^0 ({\Oc}_C(\tilde{H}))^{\vee} \otimes {\Oc}_C(\tilde{H}) 
 \to \psi^*(\T_{\PR}) \to 0,
\end{equation}gives
\begin{equation}\label{eq:3.**} 
h^0(\psi^*(\T_{\PR})) = (r+1)^2 -1 = dim (PGL(r+1 , \CC)).
\end{equation}
}
\end{remark}

\noindent
Therefore, from Remark \ref{rem:5.expmod} and from Proposition 
\ref{prop:chiave} (i), we deduce the following:

\begin{proposition}\label{prop:5.nuova}
Let $S \subset \PR$ be a smooth, non-degenerate and regular surface of general 
type and let $[X] \in V_{|D|, \delta}$ be a regular point corresponding to a 
non-degenerate and g.l.n. curve in $\PR$ for which $\delta$ is admissible and the 
Brill-Noether map $\mu_{0,C}$ of ${\Oc}_C(\psi^*(H))$ is surjective. 
Then, the morphism $\pi_{|D|, \delta}$ has injective differential at $[X]$. In particular, if 
$[X]$ is the general point of an irreducible component $V$ of $V_{|D|, \delta}$, then 
$V$ has the expected number of moduli.
\end{proposition}Our next aim is to find conditions guaranteeing that 
$X$ is g.l.n. with Brill-Noether map $\mu_{0,C}$ surjective. 
We start by considering the following crucial remark.

\begin{remark}\label{rem:0reg}
\normalfont{
Suppose that $|D|$ is a complete 
linear system on $S$ whose general element is a smooth, irreducible and 
non-degenerate curve (so that $|H-D| \neq \emptyset$). 
Assume that $[X] \in V_{| D |, \delta}$ 
corresponds to a g.l.n. curve on $S$. Denote by 
$\mu : \tilde{S} \to S$ the blow-up of $S$ along $N = Sing(X)$, 
so that $\mu |_C = \varphi$, and consider $B = \sum_{i=1}^{\delta} E_i$ 
the $\mu$-exceptional divisor.

\noindent
(a) By the hypotheses on $S$ and $X$, we have
$$H^0(\Oc_{\tilde{S}} (\mu^*(H))) \cong H^0(\Oc_{S} (H)) \cong 
H^0({\Oc}_C(\tilde{H})).$$

\noindent
(b) From the exact sequence
$$0 \to \Oc_{\tilde{S}}(K_{\tilde{S}}) \to 
\Oc_{\tilde{S}}(K_{\tilde{S}} + C) \to \omega_C \to 0,$$we get 
that $H^0(\Oc_{\tilde{S}}(K_{\tilde{S}} + C)) \to H^0(\omega_C)$ is 
surjective since, by Serre duality and by hypothesis on $S$, 
$h^1(\Oc_{\tilde{S}}(K_{\tilde{S}})) = h^1(\Oc_{\tilde{S}}) = 
h^1(\Oc_S)=0$. Therefore, by linear equivalence,
$$H^0(\Oc_{\tilde{S}}(\mu^*(K_S + D) - B)) \to H^0(\omega_C)$$is surjective.

\noindent
(c) As in (b), since $h^1(\Oc_{\tilde{S}}(K_{\tilde{S}} - 
\mu^*(H))) = h^1(\Oc_{\tilde{S}}(\mu^*(H))) = h^1(\Oc_S(H)) =0$ 
by hypothesis on $S$, we get the surjective map
$$H^0(\Oc_{\tilde{S}}(\mu^*(K_S + D - H) - B) )\to H^0(\omega_C(-\tilde{H})).$$

\noindent
Thus, we can consider the following diagram:
\begin{displaymath}
\begin{array}{ccc}
H^0(\Oc_{\tilde{S}} (\mu^*(H))) \otimes 
H^0(\Oc_{\tilde{S}}(\mu^*(K_S + D - H) - B)) & 
\stackrel{\mu_{0,\tilde{S}}}{\longrightarrow} & 
H^0(\Oc_{\tilde{S}}(\mu^*(K_S + D) - B))\\
\downarrow &  & \downarrow \\
H^0({\Oc}_C(\tilde{H})) \otimes H^0(\omega_C(-\tilde{H})) & 
\stackrel{\mu_{0,C}}{\longrightarrow}  & 
H^0(\omega_C),
\end{array}
\end{displaymath}where the vertical maps are surjective by (a), (b) and 
(c). On the other hand, we have

\begin{displaymath}
\begin{array}{ccc}
H^0(\Oc_{\tilde{S}} (\mu^*(H))) \otimes 
H^0(\Oc_{\tilde{S}}(\mu^*(K_S + D - H) - B)) & 
\stackrel{\mu_{0,\tilde{S}}}{\longrightarrow} & 
H^0(\Oc_{\tilde{S}}(\mu^*(K_S + D) - B))\\
\downarrow &  & \downarrow \\
H^0({\Oc}_S(H)) \otimes H^0(\Ii_{N/S}(K_S + D -H )) & 
\stackrel{\mu_{0,S}}{\longrightarrow}  & 
H^0(\Ii_{N/S}(K_S + D )),
\end{array}
\end{displaymath}where the vertical maps are isomorphisms. 
Thus, $\mu_{0,C}$ is surjective if $\mu_{0,S}$ is.

Recall that, if $\Ii_{N/S}(K_S + D -H )$ is a 
{\em $0$-regular} coherent sheaf on $S$, the maps
$$H^0({\Oc}_S(H)) \otimes H^0(\Ii_{N/S}(K_S + D + (\alpha-1) H)) \to 
H^0(\Ii_{N/S}(K_S + D + \alpha H))$$are surjective, for all $\alpha \geq 0$ 
(for terminology and results on $m$-regularity see, for 
example, \cite{Mu}). Therefore, the $0$-regularity of 
$\Ii_{N/S}(K_S + D -H )$ is a sufficient 
condition for the surjectivity of $\mu_{0,S}$ (and so of $\mu_{0,C}$). 
By definition, 
the given sheaf is $0$-regular iff 
\begin{equation}\label{eq:0reg}
H^1(\Ii_{N/S}(K_S + D - 2H )) = H^2(\Ii_{N/S}(K_S + D - 3 H) )= (0).
\end{equation}
}
\end{remark}
 
Our next result determines numerical conditions on the divisor class $D$ 
and an upper-bound on the number of nodes $\delta$ implying 
$(\ref{eq:0reg})$. 

\begin{theorem}\label{thm:0reg} 
Let $S \subset \PR$ be a smooth surface and let $|D |$ be a 
a complete linear system on $S$ 
whose general element is a 
smooth, irreducible divisor. Suppose that:
\begin{enumerate}
\item[i)] $(D-3H)H >0 $;
\item[ii)] $(D-4H)^2>0$ and $D(D-4H)>0$;
\item[iii)] $\nu(D,H)< D(D-4H)-4$, where $\nu(D,H)$ is the 
Hodge number of $D$ and $H$ (see Def. \ref{def:3});
\item[iv)] $\delta < \frac{D(D-4H)+\sqrt{D^2(D-4H)^2}}{8}$.      
\end{enumerate}
If $X \sim D$ is a reduced, irreducible curve with 
only $\delta$ nodes as singular points and if $N = Sing(X)$, then 
$$h^1(\Ii_{N/S}(K_S + D - 2H )) = h^2(\Ii_{N/S}(K_S + D - 3 H)) = 0;$$in 
other words, $\Ii_{N/S}(K_S + D - H)$ is $0$-regular on $S$.
\end{theorem}
\begin{proof}We start by considering the vanishing 
$h^1(\Ii_{N/S}(K_S + D - 2H ))=0$. By contradiction, assume that 
$N$ does not impose independent conditions to 
$|K_S + D - 2H|$. Let $N_0 \subset N$ be a minimal 
$0$-dimensional subscheme of $N$ for which this property holds and 
let $\delta_0= |N_0 |$. This means that $h^1(S, \Ii_{N_0}(D+K_S-2H)) 
\neq 0$ and that $N_0$ satisfies the Cayley-Bacharach condition 
(see, for example \cite{GH}). Therefore, a non-zero element of 
$H^1(\Ii_{N_0}(D+K_S- 2H))$ gives rise to a non-trivial rank 2 vector bundle 
$\E \in Ext^1(\Ii_{N_0}(D-2H), \Oc_S)$ fitting in the following exact 
sequence
\begin{equation}
0 \to \Oc_S\to \E \to \Ii_{N_0}(D-2H)\to 0, \label{eq:2.1}
\end{equation}
with 
$c_1(\E)=D-2H$ and $c_2(\E)= \delta_0 \geq 0$. Hence
\begin{equation}
c_1(\E)^2 - 4 c_2(\E)=(D-2H)^2 - 4 \delta_0. \label{eq:2.2}
\end{equation}Since $D$ is effective and irreducible with 
$D^2 > 4HD >0$, from $ii)$ it follows that $D$ is a big and 
nef divisor (see Def. \ref{def:2}). 
By applying the Index theorem to the divisor pair 
$(D, D-4H)$ and by $iv)$, we get
$$2 D(D-4H) - 8 \delta \geq D(D-4H) + \sqrt{D^2(D-4H)^2} - 8 \delta >0.$$Therefore, 
$$c_1(\E)^2 - 4 c_2(\E) \geq (D-2H)^2 - 4 \delta= D(D-4H) - 4 \delta + 
4 H^2 > 0,$$which means that 
$\E$ is Bogomolov-unstable (see Definition $\ref{def:5}$ and Remark 
$\ref{rem:a}$), hence $h^0 (\E(-M)) \neq 0$.
Twisting $(\ref{eq:2.1})$ by $\Oc_S(-M)$, we obtain
\begin{equation}
0 \to \Oc_S(-M)\to \E(-M)\to \Ii_{N_0}(D-2H-M)\to 0. \label{eq:2.5}
\end{equation}

We claim that $h^0 (\Oc_S(-M))=0$; otherwise, $-M$ would be 
an effective divisor, therefore $-MA>0$, for each ample divisor $A$. 
From $(\ref{eq:2.3})$, it follows that $c_1(\E)=M+B$, so, 
by $(\ref{eq:2.3})$ and $(\ref{eq:2.1})$, 
\begin{equation}\label{eq:2.6}
M-B = 2M-D+2H \in N(S)^+.
\end{equation}
Thus $MH > \frac{(D-2H)H}{2}$; next by $i)$ it follows that $H(D-2H)>0$, hence
$-MH<0$. 

The cohomology sequence associated to 
$(\ref{eq:2.5})$ ensures there exists a divisor 
$\Delta \sim D-2H-M $ s.t. $N_0 \subset \Delta$ and s.t. 
the irreducible nodal curve $X \sim D$, whose set of nodes 
is $N$, is not a component of $\Delta$ (otherwise, 
$-M-2H$ would be an effective 
divisor, which contradicts the non-effectiveness of $-M$).

Next, by Bezout's theorem, we get  
\begin{equation}
X \Delta = X(D-2H-M) \geq 2 \delta_0. \label{eq:2.8}
\end{equation}On the other hand, taking $M$ maximal, we may further assume 
that the general section of $\E(-M)$ vanishes in codimension $2$. Denote
by $Z$ this vanishing-locus, thus, $c_2(\E(-M))= deg(Z) \geq 0$; moreover, 
$c_2(\E(-M))= c_2(\E) + M^2+ c_1(\E)(-M) = \delta_0 + M^2 - M(D-2H)$, 
which implies 
\begin{equation}
\delta_0 \geq M(D-2H-M). \label{eq:2.9}
\end{equation}(Note that $M^2 \geq 0$ since $2M - (D-2H) \in N^+(S)$ and 
$(D-2H)$ is effective).

By applying the Index theorem to the divisor pair $(D,\; 2M-D+2H)$, we get 
\begin{equation}
D^2(2M-D+2H)^2 \leq (D(D-H)-2D(D-2H-M))^2. \label{eq:2.10}
\end{equation}
Note now that, from hypotheses $i)$ and $ii)$
it follows that 
$D(D-2H)>0$, since $D(D-4H)>0$ hence $D^2-2HD>2HD>0$. 
From $(\ref{eq:2.8})$ and from the positivity of $D(D-2H)$, 
it follows 
\begin{equation}
D(D-2H)-2D(D-2H-M) \leq D(D-2H) - 4 \delta_0.\label{eq:2.11}
\end{equation} We observe that the left side member of (\ref{eq:2.11}) 
is non-negative, since $D(D-2H)-2D(D-2H-M) = D(2M-D+2H)$, where $D$ is 
effective and, by ($\ref{eq:2.6}$), $2M-D+2H \in N(S)^+$. 
Squaring both sides
of ($\ref{eq:2.11}$), together with ($\ref{eq:2.10}$), we find
\begin{equation}
D^2(2M-D+2H)^2 \leq (D(D-2H)-4 \delta_0)^2.\label{eq:2.12}
\end{equation}

On the other hand, by $(\ref{eq:2.9})$, we get
$$
(2M-D+2H)^2= 4(M- \frac{(D-2H)}{2})^2 =(D-2H)^2-4(D-2H-M)M \geq 
(D-2H)^2- 4 \delta_0, 
$$i.e.
\begin{equation}
(2M-D+2H)^2 \geq (D-2H)^2 - 4 \delta_0. \label{eq:2.13}
\end{equation} 
Next, we define
\begin{equation}
F(\delta_0):= 4 \delta_0^2 - 4 D(D-4H) \delta_0 + (DH)^2 - 
D^2 H^2.\label{eq:2.14} 
\end{equation} Putting together $(\ref{eq:2.12})$ and $(\ref{eq:2.13})$, it
follows that $F(\delta_0) \geq 0$.
We will show that, with our numerical hypotheses, one has 
$F(\delta_0)<0$, proving the statement. 

Indeed, the discriminant of the equation $F(\delta_0)=0$ is 
$D^2(D-4H)^2$, 
and it is a positive number, since $(D-4H)^2>0$, by $ii)$, and $D^2>0$. 
We remark that $F(\delta_0)<0$ iff $\delta_0 \in 
( \alpha(D, \; H), \; \beta(D, \; H))$, where
$$ 
\alpha(D,H)=\frac{D(D-4H)- \sqrt{D^2(D-4H)^2}}{8} \;
$$
and 
$$ 
\beta(D,H)=\frac{D(D-4H)+ \sqrt{D^2(D-4H)^2}}{8};
$$
so we have to show that, 
$\delta_0 \in ( \alpha(D, \; H), \; \beta(D, \; H))$.

From $iv)$, it follows that $\delta_0< \beta(D,H)$. Note that 
$\alpha(D,H) \geq 0$: indeed, if $\alpha(D,H)<0$ then 
$D(D-4H)<\sqrt{D^2(D-4H)^2}$, which contradicts the Index Theorem, 
since $D(D-4H)>0$ and $D$ nef. Moreover, we have $\alpha(D,H) <1$: for 
simplicity, put $t = D(D-4H)$; thus $\alpha(D,H) <1$ iff
$$(*) \;\; t-8 < \sqrt{D^2(D-4H)^2} = \sqrt{t^2 - 16 ((DH)^2 - D^2 H^2)}.$$If $t-8 < 0$, $(*)$ 
trivially holds; on the other hand, if $t - 8 \geq 0$, by squaring both sides of $(*)$ we get 
$t^2 - 16 t + 64 < t^2 - 16 \nu(D,H)$ which means $\nu(D,H) < t-4 = D(D-4H) -4$, i.e. hypothesis 
$iii)$. With analogous computations we get that $\beta(D,H)>1$, which ensures there exists 
at least a positive integral value for the number of nodes.

\noindent
In conclusion, our numerical hypotheses contradict $F(\delta_0)\geq 0$, 
therefore the assumption $h^1(\Ii_N(D-2H+K_S)) \neq 0$ leads to a 
contradiction.

For what concerns the other vanishing, i.e. 
$h^2(\Ii_{N/S}(K_S + D - 3 H)) = 0$, if we consider the exact sequence 
$$0 \to \Ii_{N/S}(K_S + D - 3 H) \to \Oc_{S}(K_S + D - 3 H) 
\to \Oc_{N}(K_S + D - 3 H) \to 0, $$by Serre duality 
we get $h^2(\Ii_{N/S}(K_S + D - 3 H)) 
= h^2(\Oc_{S}(K_S + D - 3 H)) = h^0(\Oc_{S}(- D + 3 H)) = 0$ since, 
by $i)$, $3H-D$ cannot be effective.
\end{proof}

\begin{coroll}\label{cor:0reg}
If $D \sim m H$  on $S$, with $m \geq 5$, and if 
$[X] \in V_{| mH |, \delta}$ is such that 
\begin{equation}\label{eq:bound}
\delta < \frac{m (m-4)}{4} deg(S),
\end{equation}then $\Ii_{N/S}(K_S + (m -1) H)$ is $0$-regular on $S$.
\end{coroll}

We may observe that Theorem \ref{thm:0reg} also implies the geometric 
linear normality of the curve $X$. To do this, we have to 
recall the following results 
from \cite{F2}, which are a generalization of what 
Chiantini and Sernesi proved in \cite{CS} for surfaces in $\Pt$:

\begin{theorem}\label{thm:flamgln}
Let $S \subset \PR$ be a smooth, non-degenerate and linearly normal 
surface (not necessarily of general type) such that $h^1({\Oc}_S(H))$ $ = 0$ . 
Let $|D|$ be a complete linear system on $S$ whose general 
element is supposed to be smooth, irreducible and linearly normal in $\PR$. 
Then, $X$ is g.l.n. if and only if $N$ imposes independent conditions to 
the linear system $|D + K_S - H|$ 

\end{theorem}

\begin{theorem}\label{thm:gln}
Let $S \subset \PR$ be a smooth surface and let $|D|$ be a complete linear 
system, whose general element is a smooth, irreducible divisor. 
Suppose that:
\begin{enumerate}
\item[i)] $(D-H)H>0$;
\item[ii)] $(D-2H)^2>0$ and $D(D-2H)>0$;
\item[iii)] $\nu(D,H)< 4(D(D-2H)-4)$, where $\nu(D,H)$ is the Hodge number 
of $D$ and $H$;
\item[iv)] $\delta < \frac{D(D-2H)+\sqrt{D^2(D-2H)^2}}{8}$.      
\end{enumerate}
If $X \sim D$ is a reduced, irreducible curve with 
only $\delta$ nodes as singular points and if $N = Sing(X)$, then 
$h^1(\Ii_{N/S}(D + K_S -H)) = 0$ so $N$ imposes independent conditions to 
$|D + K_S -H|$. In particular, 
if $S$ is also assumed to be non-degenerate, linearly normal and 
such that $h^1(S, \; \Oc_S(H))=0$ and if the 
general element of $|D|$ is also linearly normal in $\PR$, 
then $X$ is geometrically linearly normal.
\end{theorem}

\noindent
If $D \sim mH$ then, when $m \geq 3$, all numerical 
conditions in Theorem \ref{thm:gln} hold and $iv)$ becomes 
\begin{equation}\label{eq:bogln}
\delta < \frac{m(m-2)}{4}deg(S).
\end{equation}

\begin{remark}\label{rem:raccordo}
\normalfont{
It is a straightforward computation to verify that numerical conditions 
in Theorem \ref{thm:0reg} imply the ones in Theorem \ref{thm:gln}. Thus, if $S \subset 
\PR$ is a smooth, non-degenerate, regular and linearly normal surface of general type 
such that $h^1(\Oc_S(H))=0$ and 
if $|D|$ is a complete linear system, whose general element is a 
smooth, irreducible and 
linearly normal curve satisfying numerical hypotheses in Theorem 
\ref{thm:0reg}, then $X$ is g.l.n and the map $\mu_{0,C}$ is surjective 
(see Remark \ref{rem:0reg}). 
}
\end{remark}By summarizing, we have the following result:

\begin{theorem}\label{thm:main}
Let $S \subset \PR$ be a smooth, regular, 
non-degenerate and linearly normal surface of 
general type, 
such that $h^1(\Oc_S(H))=0$. Denote by 
$|D|$ a complete linear system, whose general element is assumed to 
be a smooth, irreducible and linearly normal curve satisfying numerical 
hypotheses in Theorem \ref{thm:0reg}. Let $[X] \in V_{| D |, \delta}$ 
be a regular point of the Severi variety (in the sense of Definition \ref{def:2.reg}), 
with $\delta$ as in $iv)$ of Theorem \ref{thm:0reg}, i.e. 
$\delta < \frac{D(D-4H)+\sqrt{D^2(D-4H)^2}}{8}$. Then, the morphism 
$$\pi_{| D |, \delta} : V_{| D |, \delta} \to {\mathcal M}_g $$has 
injective differential at $[X]$. In particular, 
$\pi_{| D |, \delta} $ 
has finite fibres on each generically 
regular component of $V_{| D |, \delta}$, 
so each such component parametrizes a family 
having the expected number of moduli. 
\end{theorem}
\begin{proof}
See Remark \ref{rem:pareschi}, Remark \ref{rem:0reg}, Theorem 
\ref{thm:0reg} and Remark \ref{rem:raccordo}.
\end{proof}

\begin{coroll}\label{cor:main}
Let $S$ be as in Theorem \ref{thm:main} and let 
$D \sim m H$ on $S$, with $m \geq 5$, and assume that 
$[X] \in V_{| D |, \delta}$ is a regular 
point of the Severi variety, with $\delta$ as in (\ref{eq:bound}), i.e. 
$$\delta < \frac{m (m-4)}{4} deg(S).$$Then, the morphism 
$\pi_{| mH |, \delta} $ has injective 
differential at $[X]$. In particular, 
$\pi_{| mH |, \delta} $ 
has finite fibres on each generically 
regular component of $V_{| mH |, \delta}$, 
so each such component parametrizes a family 
having the expected number of moduli.  
\end{coroll}

\begin{remark}
\normalfont{
A particular case of the corollary above is when 
$S$ is a complete intersection in $\PR$ of type $(a_1, \ldots, a_{r-2})$; 
as already observed the upper-bound on $\delta$, ensuring that 
$X$ is g.l.n., becomes $\delta < \frac{m (m-2)}{4} deg(S)$, 
as in (\ref{eq:bogln}), 
whereas the bound on $\delta$ ensuring 
that all components of $V_{| D |, \delta}$ are regular is
\begin{equation}\label{eq:boreg}
\delta < \frac{m (m-2( (\Sigma_{i=1}^{r-2}a_i) -r-1)}{4} deg(S)
\end{equation}(see \cite{CS} and \cite{F1}). This shows that, in general, 
the strongest restriction on $\delta$ is given by asking the regularity property 
of the point $[X]$ in the sense of Severi variety theory, 
then the $0$-regularity property of the 
sheaf $\Ii_{N/S}(D+K_S-H)$ on $S$ and, finally, 
the geometric linear normality property for the curve $X$.}
\end{remark}

\begin{remark}\label{rem:example}

\normalfont{As an interesting related result, we may observe that the 
bound on $\delta$ in Theorem \ref{thm:0reg} ensuring the 
$0$-regularity of the sheaf 
$\Ii_{N/S}(K_S + D - H)$ is sharp. The following 
example was inspired by Corollary C in \cite{Tan}.}
\end{remark}

\noindent
{\bf Example}: Let $S \subset \Pt$ be a smooth sextic. 
We want to show there exist irreducible nodal curves $X$, such that 
$[X] \in V_{| 8 H | , 48 }$, for which 
$\Ii_{N/S}(K_S + D - H) = \Ii_{N/S}(9)$ is not $0$-regular. 
Since $ X \sim 8 H$, one trivially has
$$h^2(\Ii_{N/S}(K_S + D - 3H) ) = h^2( \Ii_{N/S}(7)) 
= h^2(\Oc_S(7)) = h^0 (\Oc_S(-5))=0;$$thus 
the condition of $0$-regularity fails 
as soon as $h^1(\Ii_{N/S}(8) ) \neq 0$. We will show that, for such a curve 
$X$, its set of nodes $N$ imposes one condition less to $|8 H|$ proving the sharpness 
of (\ref{eq:bound}) in Corollary \ref{cor:0reg} (observe in fact 
that $48 = \frac{6}{4}8 ( 8-4)$). 

\noindent
As a preliminary count, observe that 
the family of curves in $|8H|$ with nodes 
in $48$ given points has, at least, dimension 10. To construct an explicit example, let 
$N$ be a 0-dimensional complete intersection subscheme 
of $S$ obtained by the intersection of a general element $C_2$ of $|2 H |$ and of a 
general element $C_4$ of $|4 H |$; thus $N$ is supported on $48$ reduced points. 
By using the Koszul sequence 
of $N$ in $S$, we immediately find
$$ h^1(\Ii_{N/S}(9)) = 0 \;\; {\rm and} \;\; h^1(\Ii_{N/S}(8)) = 1.$$Observe 
that $$dim (|\Ii_{N/S}(6)|) = 43, 
\; dim (|\Ii_{N/S}(4)|) = 10, \;  
dim (|\Ii_{N/S}(2)|) = 0;$$let 
$\Gamma_4 , \; \Lambda_4 \in | \Ii_{N/S}(4) |$, 
$ \Delta_6 \in | \Ii_{N/S}(6) |$ and 
$ \Delta_2 \in | \Ii_{N/S}(2) |$ be general 
elements in such linear systems, which are smooth curves 
simply passing through $N$. Put
$$Y_1 = \Gamma_4 + \Lambda_4  \;\; {\rm and} \;\; 
Y_2 = \Delta_2 + \Delta_6;$$thus $Y_1$ and $Y_2$ are reducible 
nodal curves on $S$,  linearly equivalent to $8H$ and 
having nodes in $N$. Let 
$${\F}_{\lambda, \mu} = 
\{ \lambda Y_1 + \mu Y_2 | [\lambda , \mu ] \in {\Pp}^1 \}$$be the pencil 
of curves 
generated by $Y_1$ and $Y_2$. Its general element 
$X_{\lambda , \mu}$ is an irreducible curve 
linearly equivalent to $8H$ on $S$ passing doubly through $N$. To conclude, 
we have to show that $X_{\lambda , \mu}$ 
has only nodes in $N$. To prove this, observe that 
$$\Gamma_4 \Delta_2 = \Lambda_4 \Delta_2 = 48;$$thus, among the points 
$Y_1 Y_2 = ( \Gamma_4 + \Lambda_4 ) (\Delta_2 + \Delta_6) $, those 
which are nodes for both 
$Y_1$ and $Y_2$ are only the points of $N$. Therefore, $X_{\lambda , \mu}$ 
has only nodes in $N$.

\noindent
On the other hand, observe that such curves are geometrically linearly normal, 
since $48$ is strictly less than the bound in (\ref{eq:bogln}) 
which is $\frac{6}{4}8 (8-2) = 72$.  

\begin{remark}\label{rem:2}
\normalfont{
Note that the example above also determines non-regular 
points of the Severi variety $V_{| 8H | , 48}$. We recall that 
Chiantini and Sernesi constructed in \cite{CS} some examples of non-regular 
points of the Severi varieties 
$V_{| mH | , \frac{5}{4}m(m-4)}$, $m \geq 5$, on a 
general quintic surface $S \subset \Pt$, proving the sharpness 
of (\ref{eq:boreg}). These examples were generalized in \cite{F2} to 
Severi varieties on general canonical (i.e. $K_S \sim H$) and 
non-degenerate complete intersection surfaces in $\PR$. The key point to 
construct such examples was that on a canonical surface the condition for 
a nodal curve $X \subset S$ to be g.l.n. is equivalent to 
the fact that $[X]$ is a regular point; in 
particular, (\ref{eq:bogln}) and (\ref{eq:boreg}) coincide. 

\noindent
In the same way, when $S$ is 2-canonical (as it particularly 
happens in the example above) 
the $0$-regularity of the sheaf $\Ii_{N/S}(K_S + D -H)$ is 
equivalent to the fact that $[X]$ is a 
regular point; in particular, (\ref{eq:bound}) and (\ref{eq:boreg}) 
coincide. 
}
\end{remark}

\section{Number of moduli for families of nodal curves on 
complete intersection surfaces of general type}\label{S:6}

To complete the overview on positive answers to the 
moduli problem for divisors of the form 
$mH$ on $S \subset \PR$, the cases $1 \leq m \leq 4$, which are not 
covered by Corollary \ref{cor:main}, must be still considered. 

\noindent
From now on, we shall focus on the case of 
$S \subset \PR$ a smooth, non-degenerate complete 
intersection surface of general type; thus, 
$$K_S \sim \alpha H,$$for some positive integer $\alpha$. 

We first consider the cases $m = 3$ and $ 4$.

\begin{theorem}\label{thm:m=3,4}
Let $[X] \in V_{| mH | , \delta}$ on $S$ be a regular point, with 
$ m \geq 3$, and assume that $K_S \sim \alpha H,$ with 
$\alpha \geq 2$. If $\delta $ is as in 
(\ref{eq:bogln}), i.e.$$\delta < \frac{m(m-2)}{4}deg(S),$$then the morphism 
$$\pi_{| mH |, \delta} : V_{| mH |, \delta} \to 
{\mathcal M}_g $$has injective differential at $[X]$. In particular, 
$\pi_{| mH |, \delta} $ 
has finite fibres on each generically regular component of $V_{| mH |, \delta}$, 
so each such component parametrizes a family 
having the expected number of moduli. 

\noindent
The same conclusion holds for the family of smooth curves $V_{| mH |, 0}$ also 
with $\alpha = 1$.
\end{theorem}
\begin{proof}
By the hypothesis on $S$ and by the facts that $m \geq 3$ and $\delta$ is as 
in (\ref{eq:bogln}), we get that $X$ is g.l.n. (see Theorem \ref{thm:gln}). 
Therefore, as in (a) of Remark \ref{rem:0reg}, 
$$h^0(\Oc_{\tilde{S}} (\mu^*(H)))= h^0(\Oc_{S} (H)) = h^0({\Oc}_C(\tilde{H}))= r+1.$$By 
combining the pull-back to $\tilde{S}$ of the Euler sequence in $\PR$ and the 
exact sequence defining $C$ in $\tilde{S}$, we get the following diagram:
\begin{displaymath}
\begin{array}{ccccccc}
  & 0 &  & 0 & & 0 & \\
  & \downarrow &  & \downarrow & & \downarrow & \\
0 \to & \Oc_{\tilde{S}}(-C) & \to & H^0(\Oc_{\tilde{S}}(\mu^*(H)))^{\vee} \otimes \Oc_{\tilde{S}}(\mu^*(H)-C) 
& \to & \mu^*({\T}_{\PR}) \otimes \Oc_{\tilde{S}}(-C) & \to 0 \\
  & \downarrow &  & \downarrow & & \downarrow & \\
0 \to & \Oc_{\tilde{S}} & \to & H^0(\Oc_{\tilde{S}}(\mu^*(H)))^{\vee} \otimes \Oc_{\tilde{S}}(\mu^*(H)) 
& \to & \mu^*({\T}_{\PR}) & \to 0 \\
 & \downarrow &  & \downarrow & & \downarrow & \\
0 \to & \Oc_{C} & \to & H^0(\Oc_{C}(\tilde{H}))^{\vee} \otimes \Oc_{C}(\tilde{H}) 
& \to & \varphi^*({\T}_{\PR}) & \to 0 \\
  & \downarrow &  & \downarrow & & \downarrow & \\
  & 0 &  & 0 & & 0 & . 
\end{array}
\end{displaymath}From the regularity of $S$, we get 

\begin{displaymath}
\begin{array}{ccccl}
  & \downarrow & & \downarrow &   \\

\cdots \to & H^0(\Oc_{\tilde{S}}(\mu^*(H)))^{\vee} \otimes H^0(\Oc_{\tilde{S}}(\mu^*(H))) 
 & \stackrel{g}{\to} & H^0(\mu^*({\T}_{\PR})) & \to 0 \\

 & \downarrow^{h} & & \downarrow^{h'} &   \\

\cdots \to & H^0(\Oc_{C}(\tilde{H}))^{\vee} \otimes H^0(\Oc_{C}(\tilde{H}) )
& \stackrel{g'}{\to} & H^0(\varphi^*({\T}_{\PR})) & \to H^1(\Oc_{C}) \\

  & \downarrow & & \downarrow &  \\

  &  & & H^1( \mu^*({\T}_{\PR}) \otimes \Oc_{\tilde{S}}(-C)) & .  
\end{array}
\end{displaymath}From the second row, $\mu_{0,C}$ is surjective if and only if $g'$ is. 
Since $h'\circ g = g' \circ h $ and since $g$ is surjective, it suffices to prove that 
$h'$ is surjective. With the given hypotheses, we shall prove that 
$$(*) \;\;\;\;\;\; h^1(\mu^*({\T}_{\PR}) \otimes \Oc_{\tilde{S}}(-C)) =0$$holds. 
By Serre duality and by Leray's isomorphism, 
$ h^1(\mu^*({\T}_{\PR}) \otimes \Oc_{\tilde{S}}(-C)) = 
h^1(\Ii_{N/S} \otimes \Omega^1_{\PR}|_S \otimes \Oc_S(K_S + mH)).$ The regularity of 
$V_{| mH | , \delta} $ at $[X]$ implies that the restriction map 
$$H^0(\Oc_S(mH)) \stackrel{\rho_m}{\to} H^0(\Oc_N(mH))$$is surjective. By tensoring 
the exact sequence 
$$0 \to \Ii_{N/S} \to \Oc_S \to \Oc_N \to 0$$with the vector 
bundle $ \Omega^1_{\PR} |_S \otimes \Oc_S(K_S + mH)$, we get 
$$ \cdots \to H^0( \Omega^1_{\PR} |_S \otimes \Oc_S(K_S + mH)) 
\stackrel{\rho_{\Omega^1_{\PR} |_S \otimes \Oc_S(K_S + mH)}}{\longrightarrow} 
H^0( \Omega^1_{\PR} |_N \otimes \Oc_S(K_S + mH)) \to $$ 
$$H^1(\Ii_{N/S} \otimes \Omega^1_{\PR} |_S \otimes \Oc_S(K_S + mH)) \to 
H^1(\Omega^1_{\PR} |_S \otimes \Oc_S(K_S + mH)) \to \cdots.$$Since $S$ is a 
non-degenerate c.i. (in particular projectively normal), from standard computations 
involving the Euler sequence restricted to $S$ we find 
$ h^1(\Omega^1_{\PR} |_S \otimes \Oc_S(K_S + mH))= 0$ (for details, see \cite{F}). 
Therefore, the vanishing $(*)$ holds if and only if the map $\rho_{\Omega^1_{\PR} |_S \otimes \Oc_S(K_S + mH)}$ 
is surjective. By the assumption $K_S \sim \alpha H$, with $\alpha \geq 2$, 
the vector bundle 
$\Omega^1_{\PR} |_S \otimes \Oc_S(\alpha)$ is globally generated; 
then one concludes as in Theorem \ref{thm:5.vannodi}. In the same way, one concludes also in the case $\alpha=1$ 
and $\delta =0$.  

\noindent
Since we have proven that $X$ is g.l.n and that the map $\mu_{0,C}$, 
by Propositions \ref{prop:chiave} (i) and \ref{prop:5.nuova} we get the statement. 
\end{proof}The above result gives new positive answers to the moduli problem 
for Severi varieties of the form $V_{|mH|, \delta}$, for $m =3$ and $4$ on smooth, complete 
intersection surfaces of general type. These cases are covered neither 
by the results in Section \ref{S:3} nor by those in Section \ref{S:5}.

For what concerns the cases $m = 1$ and $2$, we cannot apply Theorem 
\ref{thm:m=3,4}, since by hypothesis, $m$ must be bigger than 2. In 
such cases, we shall make use of the following theorem in \cite{Ciro} 
(which generalizes a result of Accola in \cite{A}):

\begin{theorem}\label{thm:ciro}(see \cite{Ciro}, Teorema 2.11)
Let $\Gamma \subset \PR$ be an irreducible, non-degenerate curve of 
degree $n$ and let 
$\pi : \tilde{\Gamma} \to \Gamma $ be its normalization. 
Let $n \geq r \geq 2$ and let 
$\chi(n,r)$ be the {\em Castelnuovo number}, which is a non-negative 
integer such that 
\begin{equation}\label{eq:castelnuovo}
\frac{n-1}{r-1} - 1 \leq \chi(n,r) <  \frac{n-1}{r-1},
\end{equation}where $\chi(n,r)=0$ iff $\Gamma$ is a 
smooth, rational normal curve. Put 
\begin{equation}\label{eq:gnr}
g(n,r) = \chi(n,r) [ n -r - \frac{\chi(n,r) -1}{2} (r-1)].
\end{equation}Assume there exists on $\tilde{\Gamma}$ a linear system $g^s_m$ with 
$m \leq n$ and $s \geq r$. Then, either

\noindent
(i) $g^s_m = g^r_n$, (where $g^r_n$ is the birational linear system on 
$\tilde{\Gamma}$ related to $\pi$) 

\noindent
or

\noindent
(ii) $g(\tilde{\Gamma}) \leq \Phi(n,r):= g(n,r) - \chi(n,r) + 1 $.

\end{theorem}

\begin{remark}\label{rem:ciro}
\normalfont{In our cases, we have that $\Gamma = X $ 
is a nodal curve which is linearly equivalent 
to $mH$ on a smooth, complete intersection surface 
$S \subset \PR$ of degree $d$,  $\tilde{\Gamma} = C$ and 
$\pi = \varphi$. 

\noindent
(a) When $m=2$, $X \subset \PR$ is a 
non-degenerate, irreducible, nodal curve of degree $2d$ on $S$ and $\varphi$ is 
related to a linear system 
$g^r_{2d}$ mapping $C$ birationally onto $X$. 
By adjunction on $S$, 
$$g(C) = p_a(X) - \delta = \frac{(2H + K_S) 2H}{2} + 1 - 
\delta = (2+\alpha)d + 1 - \delta.$$From 
Theorem \ref{thm:ciro}, if $\delta < p_a(X) - \Phi(2d,r)$, i.e. 
\begin{equation}\label{eq:bociro1}
\delta <  d(2 + \alpha) + 
\frac{(r+3-4d)}{2} \chi(2d, r) + \frac{(r-1)}{2} (\chi(2d, r))^2, 
\end{equation}with $\alpha \geq 1$ and $\chi(2d, r) \in {\ZZ}_{\geq 0} \cap [\frac{2d-1}{r-1} - 1, \frac{2d-1}{r-1})$, 
then the $g^r_{2d}$ on $C$ is uniquely determined.  

\noindent
(b) If $m=1$, we have a nodal curve $X \sim H$ on $S$, which is a 
hyperplane section of a non-degenerate surface, so 
$X \subset {\Pp}^{r-1}\cong H$ is 
non-degenerate in $H$. Thus, we have a $g^{r-1}_{d}$ 
on $C$. As before, if $\delta < p_a(H) - \Phi(d,r-1) $, i.e. 
\begin{equation}\label{eq:bociro2}
 \delta < \frac{d(1+\alpha)}{2} + 
\frac{(r- 2d + 2)}{2} \chi(d, r-1) + \frac{(r-2)}{2} (\chi(d, r-1))^2, 
\end{equation}with $\alpha \geq 1$ and 
$\chi(d, r-1) \in {\ZZ}_{\geq 0} \cap [\frac{d-1}{r-2} - 1, \frac{d-1}{r-2})$, the $g^{r-1}_{d}$ on $C$ is unique.
}
\end{remark}

By using Remark \ref{rem:ciro}, we can conclude with the following

\begin{theorem}\label{thm:finale}
Let $D \sim m H$ on $S$, with $ 1 \leq m \leq 2$, and assume 
that 
$[X] \in V_{| D |, \delta}$ is a regular 
point of the Severi variety. 
Suppose that
$$\delta <  d(2 + \alpha) + 
\frac{(r+3-4d)}{2} \chi(2d, r) + \frac{(r-1)}{2} (\chi(2d, r))^2 \; {\rm and} \; 
\alpha \geq 1, \; 
{\rm if } \; \; m=2,$$where 
$\chi(2d, r) \in {\ZZ}_{\geq 0} \cap [\frac{2d-1}{r-1} - 1, \frac{2d-1}{r-1})$, 
and $$ \delta < \frac{d(1+\alpha)}{2} + 
\frac{(r- 2d + 2)}{2} \chi(d, r-1) + \frac{(r-2)}{2} (\chi(d, r-1))^2 \; {\rm and} \; 
\alpha \geq 1, 
\; {\rm if } \; m=1,$$where $\chi(d, r-1) \in 
{\ZZ}_{\geq 0} \cap [\frac{d-1}{r-2} - 1, \frac{d-1}{r-2})$. 
Then, the morphism $\pi_{| mH |, \delta} $ has injective 
differential at $[X]$. In particular, 
$\pi_{| mH |, \delta} $ 
has finite fibres on each generically regular component of $V_{| mH |, \delta}$, 
so each such component parametrizes a family 
having the expected number of moduli.  
\end{theorem}
\begin{proof}
Suppose, by contradiction, that $h^0(\varphi^*(\T_S)) \neq 0$; thus, 
$dim (\pi_{| mH |, \delta}^{-1}(\pi_{| mH |, \delta}([X]))) 
> 0$. Since $[X]$ is by assumption a regular point, it corresponds 
to an unobstructed curve in $S$. Therefore, an element of 
$T_{[X]} (\pi_{| mH |, \delta}^{-1}(\pi_{| mH |, \delta}([X]))$ is 
induced by an effective algebraic deformation. From what observed in 
Remark \ref{rem:ciro}, such deformations 
must be induced by projectivities. Then, one can conclude by using 
Proposition \ref{prop:chiave}.
\end{proof}

\noindent
{\bf Example}: if we consider an irreducible, nodal plane section $X$ on 
a smooth quintic $S \subset \Pt$, we get that $\chi(d, r-1) = 
\chi(5,2) = 3$; so 
if $[X]$ is a regular point of the 
corresponding Severi variety and, 
by (\ref{eq:bociro2}), if $\delta < \frac{10}{2} - 
\frac{15}{2} + \frac{9}{2} = 2$, the component passing through 
$[X]$ has the expected number of moduli.

\begin{remark}\label{rem:cojo}
\normalfont{
We cannot apply what observed in Remark \ref{rem:ciro} when $m=3$ and 
$4$ since, in such cases, one can show that 
$p_a(3H) - \Phi(3d,r) <0$ and $p_a(4H) - \Phi(4d, r) <0$.
}
\end{remark}

\section{Examples and final remarks}\label{S:7}

For clarity sake, here we shall summarize what one can deduce from our more 
general results of Sections 3, 5 and 6 in the particular cases of 
Severi varieties $V_{|mH|, \delta}$ on $S \subset \PR$ a smooth, 
non-degenerate complete intersection of general type or, in particular, on 
$S = S_d \subset \Pt$ of degree $d \geq 5$.

\begin{proposition}
Let $S \subset \PR$ be a smooth, non-degenerate complete intersection 
of general type whose canonical divisor is $K_S \sim \alpha H$, 
where $H$ denotes its hyperplane section. Suppose that $[X]$ is a regular 
point of the Severi variety $V_{| mH |, \delta}$.

\noindent
Assume that:  

\noindent
(1) $\delta \leq {\rm dim} (|mH|)$ if 
\begin{itemize}
\item[a)] $\alpha \geq 2$, $m \geq \alpha + 6$, $\delta \geq 1$ or
\item[b)] $\alpha \geq 1$, $m \geq \alpha + 6$, $\delta= 0$; 
\end{itemize}

\noindent
(2) $\delta < \frac{m(m-4)}{4} deg(S)$ if $\alpha \geq 1$ and $5 \leq m \leq \alpha + 5$;

\noindent
(3) \begin{itemize}
\item[a)] $\delta < \frac{m(m-2)}{4} deg(S)$ if $\alpha \geq 2$ and $m = 3, \; 4$ or
\item[b)] $\delta = 0$ if $\alpha =1$ and $m = 3, \; 4$; 
\end{itemize}

\noindent
(4) $\delta <  deg(S) (2 + \alpha) + 
\frac{(r+3-4deg(S))}{2} \chi(2deg(S), r) + \frac{(r-1)}{2} (\chi(2deg(S), r))^2$ 
if $\alpha \geq 1$ and $m =2$, where 
$\chi(2deg(S), r)$  is a non-negative integer in 
$[\frac{2deg(S)-1}{r-1}-1 , \; \frac{2deg(S)-1}{r-1})$; 

\noindent
(5) $\delta < \frac{deg(S)}{2} (1+ \alpha) + 
\frac{(r- 2deg(S) + 2)}{2} \chi(d, r-1) + \frac{(r-2)}{2} (\chi(d, r-1))^2$ if 
$\alpha \geq 1$ and $m=1$, 
where  $\chi(deg(S), r)$ is a non-negative integer in $[\frac{deg(S)-1}{r-2}-1 , 
\; \frac{deg(S)-1}{r-2})$;

\vskip 0,2cm

\noindent
Then the morphism $$\pi_{ |mH|, \delta} : V_{|mH|, \delta} 
\to {\mathcal M}_g $$has 
injective differential at $[X]$. In particular, it has 
finite fibres on each generically regular component of $V_{|mH|, \delta}$, 
so each such component parametrizes a family 
having the expected number of moduli.
\end{proposition} 

\vskip 0,2cm

In particular, we have:

\begin{coroll}
Let $S \subset \Pt$ be a smooth surface of degree $d \geq 5$ and let 
$[X] \in V_{| mH |, \delta}$ be a regular point.

\noindent
Assume that:  

\noindent
(1) $\delta \leq {\rm dim} (|mH|)$ if 
\begin{itemize}
\item[a)] $d \geq 6$, $m \geq d+2$, $\delta \geq 1$ or
\item[b)] $d \geq 5$, $m \geq d+2$, $\delta= 0$; 
\end{itemize}

\noindent
(2) $\delta < \frac{m(m-4)}{4} d$ if $d \geq 5 $ and $5 \leq m \leq d+1$;

\noindent
(3) \begin{itemize}
\item[a)] $\delta < \frac{m(m-2)}{4} d$ if $d \geq 6 $ and $m = 3, \; 4$ or
\item[b)] $\delta = 0$ if $d = 5 $ and $m = 3, \; 4$; 
\end{itemize}

\noindent
(4) $\delta <  d - 2 $ if $d \geq 5$ and $m =2$; 

\noindent
(5) $\delta < d - 3 $ if $ d \geq 5$ and $m=1$.

\vskip 0,2cm

\noindent
Then the morphism $$\pi_{ |mH|, \delta} : V_{|mH|, \delta} 
\to {\mathcal M}_g $$has 
injective differential at $[X]$. In particular, it has 
finite fibres on each generically regular component of $V_{|mH|, \delta}$, 
so each such component parametrizes a family 
having the expected number of moduli.
\end{coroll}

Observe that our results generalize what can be proven 
in the case of a general smooth, 
complete intersection surface $S \subset \PR$ 
by using a recent result of Schoen, \cite{Schoen}. In his paper, he studies 
algebraic varieties which are dominated by products 
of varieties of smaller dimension (abbreviated {\em DPV}); in the case of products of curves, 
one writes {\em DPC}. The main goal of Schoen's paper is to discuss, via 
real algebraic group theory and Hodge theory, some 
obstructions to {\em DPC} and {\em DPV} properties. As a result, he shows 
for example that if $W \subset {\Pp}^N$ is a sufficiently general complete 
intersection variety of degree $d > N+1$ and of dimension $n \geq 2$, then $W$ cannot satisfy 
the {\em DPC}-property. 
Thus, the general complete intersection surface $S \subset \PR$, of degree $d \geq r+2$, 
cannot be dominated by a product of curves $C_1 \times C_2$. Therefore, 
there cannot exist isotrivial pencils of smooth or $\delta$-nodal 
curves in $|mH|$, otherwise, after a suitable base change, 
such a surface would be {\em DPC}.

Thus, via Schoen's results, one can 
answer the moduli problem, for smooth and nodal 
curves in the linear system $|m H|$, $m \geq 1$, on a 
general complete intersection surface $S \subset \PR$ of degree 
$d \geq r+2$. Our results are more generally 
valid for divisors $D$ on $S$, where $S$ is not necessarily a 
general complete intersection and can have a wildly complicated 
$Div(S)$.

\end{document}